\documentclass[12pt,twoside,a4paper]{amsart}
\UseRawInputEncoding
\usepackage{amsmath,amssymb}
\usepackage{cite}
\usepackage{graphicx}
\usepackage{txfonts}
\usepackage{mathrsfs}
\usepackage{empheq}
\newcommand\res{\mathop{\mathrm{res}}}
\newcommand{\ds}{\displaystyle}
\newcommand{\K}{J.Kaczorowski and K.Wiertelak}
\newcommand{\ka}{J.Kaczorowski}

\newcommand{\s}{\mathcal{S}^\text{poly}}
\newcommand{\se}{\mathcal{S}}
\newcommand{\f}{\overline{F(1-\overline{s})}}
\newcommand{\mf}{\overline{\mu_F (k)}}

\allowdisplaybreaks[4]
\theoremstyle{definition}
\newtheorem{definition}{Definition}[subsection]
\newtheorem{remark}[definition]{Remark}

\theoremstyle{theorem}
\newtheorem{theorem}{Theorem}[subsection]
\newtheorem{lemma}[theorem]{Lemma}

\newtheorem{fact}[theorem]{Fact}
\newtheorem{cor}[theorem]{Corollary}
\begin{document}

\begin{center}
\textbf{\large  An analytic approach to the remainder terms in the asymptotic formulas \\ - the Volterra integral equation, the Whittaker function - }
\\
\vspace{0.9cm}
by
\\
\vspace{0.35cm}
Hideto  IWATA
\end{center}

\vspace{0.5cm}

\textbf{Abstract.} In the present paper, firstly, we consider the Volterra integral equation of second type for a remainder term in an asymptotic formula of an arithmetic function which satisfies some special conditions and obtained a solution of the equation. The method using there is applied to the remainder term in an asymptotic formula of the associated Euler totient function. Secondly, we consider a function defined a series over all non-trivial zeros of the generalized $L$-functions. We proved some analytic properties which satisfy some conditions. In particular, we use the Whittaker function which is kind of the confluent hypergeometric function there.

\section{Introduction}
\subsection{Previous research for the Euler totient function}

\hspace{1cm}

For $n \in \mathbb{N}$, let $\varphi(n)$ be the number of positive integers not exceeding $n$ which are relatively prime to $n$. The function $\varphi(n)$ is called the \textit{Euler totient function} and appears in various fields e.g. elementary number theory, group theory. In analytic number theory, studying the arithmetic mean of $\varphi(n)$ is a classical problem. Let
\begin{equation}
   E(x) = \sum_{n \leq x} \varphi(n) - \frac{3}{\pi^2}x^2   \tag{1.1.1}
\end{equation}be the associated error term. The error term (1.1.1) has been studies for a long time. P. G. Dirichlet proved the estimate
\begin{equation*} 
   E(x) \ll x^{1+\epsilon}   \tag{1.1.2}
\end{equation*}for every positive $\epsilon$. Here, we use the notation $f(x) \ll g(x)$, if there is a constant $C>0$ such that $|f(x)| \leq Cg(x)$ for all $x$ in the appropriate range. If the implicit constant may depend on $\epsilon$ in $f(x) \ll g(x)$, we use the notation $f(x) \ll_\epsilon g(x)$. Also, we use the notation $f(x) \gg g(x)$, if there is a positive constant $c$ such that $f(x) \geq cg(x)$ and $g$ is non-negative. The estimate (1.1.2) was improved by F. Mertens to $E(x) \ll x\log x$ (see~\cite{me}). A. Walfisz proved 
\begin{equation*}
E(x) \ll x(\log x)^\frac{2}{3}(\log\log x)^\frac{4}{3}   \tag{1.1.3}
\end{equation*}(see ~\cite{W}). 

\footnote[0]
{
\textit{2010 Mathematics \hspace{-0.1cm} Subject Classification Primary} : 45D05 ; 11A25 ; 11N37 ; 33C15, 11M06, 11M41 
\\
\textit{Key words and phrases}
 : Volterra integral equation of second type, the remainder term in the asymptotic formula, arithmetical function, twisted Euler $\varphi$-function, associated Euler totient function, polynomial Euler product, Whittaker function, Selberg class.
}The estimate (1.1.3) is the best known result. S.S. Pillai and S.D. Chowla proved  
\[ {\ds} \sum_{n \leq x} E(n) \sim \frac{3}{2\pi^2}x^2 \](see ~\cite{pc}) and S.D. Chowla also proved 
\[ \int_{0}^x E(t)^2 dt \sim \frac{x^3}{6\pi^2} \](see~\cite{c}). H.L. Montgomary proved the best $\Omega$-result for (1.1.1)
\begin{equation*}
   E(x) = \Omega_{\pm}(x\sqrt{\log\log x})  \tag{1.1.4}
\end{equation*}(see~\cite{Mon}). Here, we use the notation $f(x) = \Omega_+ (g(x))$ and  $f(x) = \Omega_- (g(x))$ i.e. the inequalities $f(x) > Cg(x)$ and $f(x) < -Cg(x)$ hold respectively for some arbitrarily large values of $x$ and a suitable positive constant $C$. Also, we use the notation $f(x) = \Omega_\pm (g(x))$ i.e. both $f(x) = \Omega_+ (g(x))$ and $f(x) = \Omega_- (g(x))$ hold. In 2010, {\K} considered
\begin{equation*}
   \tilde{E}(x) := \int_{0}^x E(t)\frac{dt}{t} = \sum_{n \leq x}\varphi(n)\log\frac{x}{n} - \frac{3}{2\pi^2}x^2.  \tag{1.1.5}
\end{equation*} The right hand side in (1.1.5) is obtained by Fact 2.3.4 in Section 2.3. {\K} proved that for $x$ tending to infinity  
\begin{equation*}
    \tilde{E}(x) = \Omega_{\pm} (x^\frac{1}{2}\log\log\log x).
\end{equation*}Moreover, under the Riemann Hypothesis for the Riemann zeta function $\zeta(s)$ i.e. $\zeta(s) \neq 0 \> \left(s=\sigma+it, \> \left( \sigma > \frac{1}{2} \right) \right)$, there exist a positive constant $B$ such that
\begin{equation*}
    \tilde{E}(x) \ll x^\frac{1}{2}\exp\left( B\frac{\log x}{\log\log x} \right)
\end{equation*}(see~\cite{K}). We write a complex variable $s=\sigma + it$ in what follows. We use the same notation throughout this paper. {\K} also studied (1.1.1) by splitting into two summands. {\K} considered there the following Volterra integral equation of second type for (1.1.1) (see~\cite{Kac and Wie}) :
\begin{equation}
   F(x) - \int_{0}^\infty K(x,t)F(t)dt = E(x)  \quad (x \geq 1), \tag{1.1.6}
\end{equation}where $F(x)$ is the unknown function and the kernel $K(x,t)$ is defined as follows:
\begin{equation*}
   K(x,t) = \begin{cases}
                 1/t \quad (0 < t \leq x),
                 \\
                  0  \quad (1 \leq x < t).   \tag{1.1.7} 
              \end{cases}
\end{equation*}The equation (1.1.6) can be solved explicitly. Let us put
\begin{equation*}
   f(x) = -\sum_{n=1}^\infty \frac{\mu(n)}{n}\left\{ \frac{x}{n} \right\}   \tag{1.1.8}
\end{equation*}for every $x \geq 0$, where $\mu(n)$ denotes the M\"{o}bius function i.e.
\begin{equation*}
   \mu(n) = \begin{cases}
                      1 \hspace{1cm} (n=1),
                      \\
                      (-1)^r \hspace{0.3cm} (n=p_1 p_2 \cdots p_r, p_i  \> (1 \leq i \leq r) : \text{prime}),
                      \\
                      0 \hspace{1cm} (\text{otherwise}),
                  \end{cases}
\end{equation*}and $\{x\}=x-[x]$ is the fractional part of a real number $x$. 
\begin{theorem}[Theorem 1.1 in ~\cite{ Kac and Wie }]
The general solution of (1.1.6) is
\begin{equation*}
   F(x) = (f(x)+A)x,   \tag{1.1.9}
\end{equation*}where $A$ is an arbitrary constant. 
\end{theorem}In~\cite{Kac and Wie}, $F(x)=xf(x)$ is claimed to be the unique solution of the integral equation (1.1.6), but this uniqueness does not hold even assuming the initial value condition at $x=0$. Probably, the term $Ax$ is missing to give the general solution. {\K} defined the linear space $\mathcal{X}$ as follows :
\footnotesize
\[  \mathcal{X} = \left\{ g : (0,\infty) \to \mathbb{R} ;  \text{Lebesgue locally integrable},
     \int_{0}^1 |g(t)||\log t|^N \frac{dt}{t} < \infty \quad (\forall N \in \mathbb{N}) \right\}. \]
\normalsize
Also, {\K} defined the operator $\delta_1$ on $\mathcal{X}$ as follows :
\begin{equation*}
   \delta_1 (g)(x) = \int_{0}^x \frac{g(t)}{t}dt \quad (x>0, g \in \mathcal{X}).   \tag{1.1.10}
\end{equation*}Moreover, let $\delta_k $ denote the $k$-fold iteration of $\delta_1$ : 
\begin{equation*}
   \delta_k = \delta_1 \circ \cdots \circ \delta_1 \quad (k \> \text{times})   \tag{1.1.11}
\end{equation*}and let
\begin{equation*}
   R_k (x) = \delta_k (E)(x),   \tag{1.1.12}
\end{equation*}where $E(x)$ is defined in (1.1.1). For $x \geq 0$ let us write 
\begin{equation*}
   g(x) = \sum_{n=1}^\infty \mu(n)\left\{ \frac{x}{n} \right\}^2.  \tag{1.1.13}
\end{equation*}
\begin{theorem}[Theorem 1.2 in ~\cite{ Kac and Wie }]
For $x \geq 1$ we have
\[  E(x) = xf(x) + \frac{1}{2}g(x) + \frac{1}{2}. \]
\end{theorem}
According to Theorem1.1.2, for $x \geq 1$ we can split $E(x)$ as follows : 
\begin{equation*}
   E(x) = E^{\text{AR}}(x) + E^{\text{AN}}(x),  \tag{1.1.14}
\end{equation*}where
\begin{equation*}
   E^{\text{AR}}(x) = xf(x), \quad \text{and} \quad E^{\text{AN}}(x) = \frac{1}{2}g(x) + \frac{1}{2}
\end{equation*}with $f(x)$ and $g(x)$ given by (1.1.8) and (1.1.13) respectively. We call $E^{\text{AR}}(x)$ and $E^{\text{AN}}(x)$ the \textit{arithmetic} and the \textit{analytic part} of $E(x)$ respectively. J.Kaczorowski and K.Wiertelak proved the $\Omega$-estimates for $E^{\text{AR}}(x)$ and $E^{\text{AN}}(x)$. On an arithmetic part, {\K} obtained an $\Omega$-result for a class of arithmetic functions $\mathcal{A}$ as follows (see ~\cite{ Kac and Wie }) : let $\mathcal{A}$ denote the set of all arithmetic functions $\alpha(n)$ satisfying the following conditions :  
\begin{enumerate}
\item[(i)] $\alpha(n)$ is real and multiplicative.
\item[(ii)] There exists a positive real number $\theta < 1$ such that
               \begin{equation*}
                   \alpha(n) \ll n^\theta .
                \end{equation*}
\item[(iii)] We have
               \begin{equation*}
                  \sum_{n=1}^\infty \frac{\alpha(n)}{n^2} \neq 0.  
               \end{equation*}
\item[(iv)] For every $N \geq 1$ we have
                \begin{equation*}
                   \sum_{n=1}^N |\alpha(n)| \ll N.
                \end{equation*}
\item[(v)] The series
                \begin{equation*}
                   \sum_{n=1}^\infty \frac{\alpha(n)}{n}
                \end{equation*}converges.
\item[(vi)] There exists a positive real number $\eta$ and a sequence of positive numbers $x_\nu \to \infty$ such that
                \begin{equation*}
                   \sum_{\substack{p \leq x_\nu \\ \alpha(p) < 0 \\ p \equiv 3 \hspace{-0.1cm} \pmod{4}}} \frac{|\alpha(p)|}{p} \geq \eta\log\log x_\nu + O(1)
                 \end{equation*}for all $\nu \geq 1$.
\end{enumerate}For $\alpha(n) \in \mathcal{A}$, we write
\begin{equation*}
   f(x,\alpha) = \sum_{d=1}^\infty \frac{\alpha(d)}{d}s\left( \frac{x}{d} \right),
\end{equation*}where $s(x)$ denotes the saw tooth function :
\begin{equation*}
   s(x) = \begin{cases}
                 0 \quad (x \in \mathbb{Z}),
                 \\
                 \frac{1}{2} - \{ x \} \quad (\text{otherwise}).  \tag{1.1.15}
             \end{cases}
\end{equation*}Since
\begin{equation*}
   \sum_{n=1}^\infty \frac{\mu(n)}{n} = 0,
\end{equation*}we have
\begin{equation*}
   f(x,\mu) = \frac{1}{2}(f(x-0) + f(x+0)),
\end{equation*}where $f(x)$ is the same as in (1.1.8). Hence, for $x \notin \mathbb{Z}$ we have $f(x,\mu) = f(x)$. Moreover, for every $x \geq 1$ let us put
\begin{equation*} 
   R(x,\alpha) = \sup_{y \geq x} \left| \sum_{n>y} \frac{\alpha(n)}{n} \right|
\end{equation*}and
\begin{equation*}
   R^\ast (x,\alpha) = \sqrt{R(\sqrt{x},\alpha)} + \frac{1}{x}.
\end{equation*}
\begin{theorem}[Theorem 1.3 in ~\cite{ Kac and Wie }]
Let $\alpha \in \mathcal{A}$. Then we have
\begin{equation*}
   f(x,\alpha) = \Omega_\pm \left( \left( \log\log\frac{1}{R^\ast (x^\frac{2}{3},\alpha)} \right)^\eta \right)
\end{equation*}as $x \to \infty$.
\end{theorem}
\begin{cor}[Corollary 1.4 in ~\cite{Kac and Wie}]
Let $f(x)$ be defined in (1.1.8). Then
\begin{equation*}
   f(x) = \Omega_{\pm} (\sqrt{\log\log x}) \quad \text{and} \quad E^{\text{AR}}(x) = \Omega_{\pm} (x\sqrt{\log\log x})   \tag{1.1.16}
\end{equation*}as $x \to \infty$.
\end{cor}On an analytic part, {\K} obtained an $\Omega$-result as follows :
\begin{theorem}[Theorem 1.8 in ~\cite{ Kac and Wie }]
For $x$ tending to infinity we have
\begin{equation*}  
   E^{\text{AN}}(x) = \Omega_{\pm}(x^\frac{1}{2} \log\log\log x).   \tag{1.1.17}
\end{equation*}
\end{theorem}The following Theorem 1.1.6 and Lemma 1.2.3 below show that while $E^{\text{AN}}(x)$ depends on the non-trivial zeros of the Riemann zeta function, $E^{\text{AR}}(x)$ does not. In particular, Theorem1.1.6 is related to a certain equivalence condition of the Riemann Hypothesis for $\zeta(s)$.
 \begin{theorem}[Theorem 1.7 in ~\cite{ Kac and Wie }]
The following statements are equivalent.
\begin{enumerate}
\item[(1)] The Riemann Hypothesis is true.
\item[(2)] There exists a positive constant $A$ such that for $x \geq e^e$ we have
\begin{equation*}  
   E^{\text{AN}}(x) \ll x^\frac{1}{2}\exp \left( A\frac{\log x}{\log\log x} \right).   \tag{1.1.18}
\end{equation*}
\item[(3)] For every $\epsilon>0$ and $x \geq 1$ we have
\[ E^{\text{AN}}(x) \ll_\epsilon x^{\frac{1}{2}+\epsilon}. \]
\end{enumerate}
\end{theorem}{\K} obtained a better decomposition for the remainder term in the asymptotic formula for a generalization of the Euler totient function ({see ~\cite{I}, ~\cite{kac and wie}) : For a non-principal real Dirichlet character $\chi \hspace{0.1cm} (\text{mod }q), q>2,\text{ let } \varphi(n,\chi)$ denote the \textit{twisted Euler $\varphi$-function} 
\begin{equation*}
\varphi(n,\chi) = n\prod_{p|n} \left( 1 - \frac{\chi(p)}{p} \right),   \tag{1.1.19}
\end{equation*}where the product in (1.1.19) is over the prime $p$ which divisors of $n$. Refer to the paper ~\cite{h} on the application of $\varphi(n,\chi)$. {\K} made a consideration similar to~\cite{Kac and Wie} for the remainder term in the asymptotic formula of the above twisted Euler $\varphi$-function. Let
\begin{equation*}
   E(x,\chi) = \sum_{n \leq x} \varphi(n,\chi) - \frac{x^2}{2L(2,\chi)}   \tag{1.1.20}
\end{equation*}and
\begin{equation*}
   E_1 (x,\chi) = \begin{cases}
                         E(x,\chi) \quad (x \notin \mathbb{N}),
                         \\
                         \frac{1}{2}(E(x-0,\chi)+E(x+0,\chi)) \quad (\text{otherwise})   \tag{1.1.21}
                      \end{cases}
\end{equation*}be the corresponding error terms. Here, as usual, $L(s,\chi)$ denotes the Dirichlet $L$-function associated to $\chi$. It is easy to see that $E(x,\chi) = O(x\log x)$ for $x \geq 2$. Let $s(x)$ be the same as in (1.1.15). We write for $x \geq 0$
\begin{equation*}
   f(x,\chi) = \sum_{d=1}^\infty \frac{\mu(d)\chi(d)}{d}s\left( \frac{x}{d} \right)   \tag{1.1.22}
\end{equation*}and
\begin{equation*}
   g(x,\chi) = \sum_{d=1}^\infty \mu(d)\chi(d) \left\{ \frac{x}{d} \right\}\left( \left\{ \frac{x}{d} \right\} - 1 \right).   \tag{1.1.23}
\end{equation*}
\begin{theorem}[Theorem1.1. in ~\cite{kac and wie}]
The solution of the following Volterra integral equation of second type
\begin{equation*}
F(x,\chi) - \int_{0}^\infty K(x,t)F(t,\chi)dt = E_1 (x,\chi) \quad (x \geq 0),   \tag{1.1.24}
\end{equation*}where $K(x,t)$ is the same as in (1.1.7) is the function
\begin{equation*}
   F(x,\chi) = (f(x,\chi)+A)x,   \tag{1.1.25}
\end{equation*}where $A$ is an arbitrary constant. 
\end{theorem}(In~\cite{kac and wie}, the unique solution is $F(x,\chi) = xf(x,\chi)$, but the comments just after (1.1.9) should also be applied here). {\K} also obtained the arithmetic and the analytic part of $E(x,\chi)$ respectively.
\begin{theorem}[Theorem1.2. in ~\cite{kac and wie}]
For $x \geq 0$ 
\begin{equation*}
   E_1 (x,\chi) = E^{\text{AR}}(x,\chi) + E^{\text{AN}}(x,\chi),  \tag{1.1.26}
\end{equation*}where
\begin{equation*}
   E^{\text{AR}}(x,\chi) = xf(x,\chi) \quad \text{and} \quad E^{\text{AN}}(x,\chi) = \frac{1}{2}g(x,\chi)   \tag{1.1.27}
\end{equation*}with $f(x,\chi)$ and $g(x,\chi)$ given by (1.1.22) and (1.1.23) respectively. 
\end{theorem}{\K} proved the $\Omega$-estimates for $E^{\text{AR}}(x,\chi)$ and $E^{\text{AN}}(x,\chi)$. Also, {\K} proved the equivalence of the Riemann Hypothesis for the Dirichlet-$L$ function in terms of $E^{\text{AN}}(x,\chi)$.
\begin{theorem}[Theorem 1.3 in ~\cite{kac and wie}]
Let $f(x,\chi)$ be defined in (1.1.22). Then
\begin{equation*}
   f(x,\chi) = \Omega_{\pm} ((\log\log x)^\frac{1}{4}) \quad \text{and} \quad E^{\text{AR}}(x,\chi) = \Omega_{\pm} (x(\log\log x)^\frac{1}{4})   \tag{1.1.28}
\end{equation*}as $x \to \infty$.
\end{theorem}

\begin{theorem}[Corollary 1.4 in ~\cite{kac and wie}]
We have
\begin{equation*}  
   E^{\text{AN}}(x,\chi) = \Omega_{\pm}( x(\log\log x)^\frac{1}{4} )    \tag{1.1.29}
\end{equation*}as $x \to \infty$.
\end{theorem}

\begin{theorem}[Theorem 1.5 in ~\cite{kac and wie}]
The following statements are equivalent.
\begin{enumerate}
\item[(1)] $L(s,\chi) \neq 0 \hspace{0.1cm} \text{for} \hspace{0.1cm} \sigma>\frac{1}{2}$.
\item[(2)] There exists a positive constant $A$ such that for $x \geq e^e$ we have
\begin{equation*}  
   E^{\text{AN}}(x,\chi) \ll x^\frac{1}{2}\exp \left( A\frac{\log x}{\log\log x} \right).   \tag{1.1.30}
\end{equation*}
\item[(3)] For every positive $\epsilon$ and $x \geq 1$ we have
\[ E^{\text{AN}}(x,\chi) \ll_\epsilon x^{\frac{1}{2}+\epsilon}. \]
\end{enumerate}
\end{theorem}On the $E^{\text{AN}}(x,\chi)$, {\K} proved the $\Omega$-estimates the cases the character is even and odd respectively (see~\cite{kac and wie}).

\subsection{Preliminaries to prove Theorem1.1.5}

\hspace{1cm}

In this section, we state lemmata to prove Theorem1.1.5.
\begin{lemma}[Lemma 5.2 in ~\cite{Kac and Wie}]
For $\sigma>2$ we have
   \begin{equation*}
      \int_{1}^\infty E^{\text{AN}}(x) x^{-s-1}dx = \frac{3}{\pi^2}\frac{1}{s-2} + \frac{\zeta(s-1)}{s(1-s)\zeta(s)}.  \tag{1.2.1}
   \end{equation*}
\end{lemma}
\begin{lemma}[Lemma 5.3 in ~\cite{Kac and Wie}]
Suppose that a measurable locally bounded function $h : [1,\infty) \to \mathbb{R}$ satisfies $h(x) = O(x^A)$, and $h(x) \leq Bx^a \log\log x$ or $h(x) \geq -Bx^a \log\log x$ for certain positive $a, A$ and $B$ and all large $x$. Moreover, let its Mellin transform $F(s) = \int_{1}^\infty h(x)x^{-s-1}dx$ be holomorphic on the interval $[a,A]$. Then $F(s)$ is holomorphic for $\sigma>a$ and 
\[ F(s) \ll \frac{1}{\sigma-a}\log\left( \frac{1}{\sigma-a} \right) \] uniformly for $a<\sigma<a+\frac{1}{2}$.
\end{lemma}We can prove the following lemma applying Lemma1.2.1 and Lemma1.2.2 to $h(x) = E^{\text{AN}}(x), a = 1/2, A=2$ :
\begin{lemma}[Lemma 5.4 in ~\cite{Kac and Wie}]
Suppose that 
\begin{equation*}
   \exists x_0 , C_0 > 0 \quad \forall x \geq x_0 \quad E^{\text{AN}}(x) \leq  C_0 x^\frac{1}{2}\log\log x  \tag{1.2.2}
\end{equation*}or
\begin{equation*}
   \exists x_0 , C_0 > 0 \quad \forall x \geq x_0 \quad E^{\text{AN}}(x) \geq  -C_0 x^\frac{1}{2}\log\log x.  \tag{1.2.3}
\end{equation*}Then the Riemann Hypothesis is true, all non-trivial zeros of the Riemann zeta function are simple, and denoting by $\rho = \frac{1}{2}+i\gamma$ a generic non-trivial zero we have
\begin{equation*}
   \frac{\zeta(\rho-1)}{\zeta^\prime (\rho)} \ll \gamma^2 \log |\gamma|.  \tag{1.2.4}
\end{equation*}
\end{lemma}Under the Riemann Hypothesis, the Lindel\"{o}f Hypothesis is true. Therefore, we have 
\begin{equation*}
   \zeta^{(k)} (s) \ll |t|^\epsilon   \tag{1.2.5}
\end{equation*}for every $\sigma \geq \frac{1}{2}$ and every integer $k \geq 0$. Using (1.2.5), the estimate (1.2.4) is proved.  
\begin{lemma}[Lemma 5.1 in ~\cite{Kac and Wie}]
With the notation (1.1.12) we have
\[  R_1 (x) + R_2 (x) = \Omega_{\pm}(\sqrt{x}\log\log\log x) \]as $x \to \infty$.    
\end{lemma}
In ~\cite{Kac and Wie}, the proof of Lemma1.2.4 is written as follows : Analogous result for $R_1 (x)$ in the place of $R_1 (x)+R_2 (x)$ was established in ~\cite{K}. The present lemma follows by repeating all steps in the proof of Theorem 1.1 in ~\cite{K}. The required modifications are straightforward and shall not be described here in detail.

\subsection{The outline of the proof of Theorem1.1.5}

\hspace{1cm}

We state the outline of the proof of Theorem1.1.5 (see ~\cite{Kac and Wie}) : We can assume that (1.2.2) and (1.2.3) are true since otherwise there is nothing left to be proved. Taking the inverse Mellin transform in (1.2.1) we obtain
\begin{equation*}
   E^{\text{AN}}(x) = \frac{1}{2\pi i}\int_{\mathcal{L}} \frac{\zeta(s-1)}{\zeta(s)}\frac{x^s}{s(1-s)}ds,   \tag{1.3.1}
\end{equation*}where the path of integration $\mathcal{L}$ consists of the half-line $[3-i\infty, 3-2i]$, the semi-circle $s=3+2e^{i\theta}$, $\frac{\pi}{2} < \theta < \frac{3}{2}\pi$ and the half-line $[3+2i, 3+i\infty]$. On the Mellin inversion, the following fact is known :
\begin{fact}[Theorem 28 in ~\cite{T}]
Let $f(y)y^{k-1} \hspace{0.01cm} (k>0)$ belongs to $L(0,\infty)$, and let $f(y)$ be of bounded variation in the neighborhood at the point $y=x$. Let 
   \begin{align*}
      \mathcal{F}(s) = \int_{0}^\infty f(x)x^{s-1}dx \hspace{0.3cm} (s=k+it). \tag{1.3.2}   
   \end{align*}Then
   \begin{align*}
      \frac{1}{2}\{ f(x+0) + f(x-0) \} = \frac{1}{2\pi i}\lim_{T \to \infty} \int_{k-iT}^{k+iT} \mathcal{F}(s)x^{-s}ds.  \tag{1.3.3}
   \end{align*}
\end{fact}
According to 
\[ \frac{1}{s(1-s)} = -\frac{1}{s^2} - \frac{1}{s^3} + \frac{1}{s^3 (1-s)}, \]we split the integral on the right hand side of (1.3.1) into three parts
\begin{align*} 
   \frac{1}{2\pi i}\int_{\mathcal{L}} \frac{\zeta(s-1)}{\zeta(s)}\frac{x^s}{s(1-s)}ds
   &= -\frac{1}{2\pi i}\int_{\mathcal{L}} \frac{\zeta(s-1)}{\zeta(s)}\frac{x^s}{s^2}ds -\frac{1}{2\pi i}\int_{\mathcal{L}} \frac{\zeta(s-1)}{\zeta(s)}\frac{x^s}{s^3}ds 
   \\
   &+ \frac{1}{2\pi i}\int_{\mathcal{L}} \frac{\zeta(s-1)}{\zeta(s)}\frac{x^s}{s^3 (1-s)}ds 
   \\
   &= -R_1 (x) -R_2 (x) + \frac{1}{2\pi i}\int_{\mathcal{L}} \frac{\zeta(s-1)}{\zeta(s)}\frac{x^s}{s^3 (1-s)}ds 
    \\
    &= -R_1 (x) -R_2 (x) + I,
\end{align*}say where $R_1 (x)$ and $R_2 (x)$ are the cases $k=1,2$ in (1.1.12) respectively. Shifting the line of integration to the left we have
\begin{equation*}
 I = \sum_{\rho} \frac{\zeta(\rho-1)}{\rho^3 (1-\rho)\zeta^\prime (\rho)}x^\rho + \frac{1}{2\pi i}\int_{\frac{1}{4}-i\infty}^{\frac{1}{4}+i\infty} \frac{\zeta(s-1)}{\zeta(s)} \frac{x^s}{s^3(1-s)}ds.   \tag{1.3.4}
\end{equation*}Using the estimate (1.2.4), for every $\epsilon>0$
\begin{align*}
   I
   &= \sum_{\rho} \frac{\zeta(\rho-1)}{\rho^3 (1-\rho)\zeta^\prime (\rho)}x^\rho + \frac{1}{2\pi i}\int_{\frac{1}{4}-i\infty}^{\frac{1}{4}+i\infty} \frac{\zeta(s-1)}{\zeta(s)} \frac{x^s}{s^3(1-s)}ds
   \\
   &\ll x^\frac{1}{2}\sum_\rho \frac{\log|\gamma|}{\gamma^2} + x^\frac{1}{4}\int_{\frac{1}{4}-i\infty}^{\frac{1}{4}+i\infty}|s|^{-3+\epsilon}|ds|
   \\
   &\ll x^\frac{1}{2} + x^\frac{1}{4}
   \\
   &\ll x^\frac{1}{2}.
\end{align*}Hence
\[ E^{\text{AN}}(x) =  -R_1 (x) -R_2 (x) + O(x^\frac{1}{2}).  \]The assertion now follows from Lemma 1.2.4. $\square$

\subsection{The analytic property of a function $f(z)$}

\hspace{1cm}

The proof of Lemma1.2.4 is just as in the case of $R_1 (x)$ in ~\cite{K} stated in Section1.2.When proving Theorem 1.1 in ~\cite{K}, {\K} used the functional equation (1.4.5) below (see p1642-3 in ~\cite{K}) : We describe basic analytic properties of the function $f(z)$ defined for Im $z >0$ as follows :
\begin{equation*}
   f(z) = \lim_{n \to \infty} \sum_{\substack{\rho \\ 0< \text{Im} \> \rho < T_n}} \frac{e^{\rho z}\zeta(\rho-1)}{\zeta^\prime (\rho)},   \tag{1.4.1}
\end{equation*}where $T_n$ denotes a sequence of real numbers yields appropriate grouping of the zeros. The summation is over non-trivial zeros the Riemann zeta-function with positive imaginary part. For simplicity we assume here that the zeros are simple. Let us denote by $\ell(-\frac{1}{4},\frac{5}{2})$ a simple and smooth curve $\tau : [0,1] \longrightarrow \mathbb{C}$ such that $\tau(0) = -\frac{1}{4}, \> \tau(1) = \frac{5}{2}$ and $0<\text{Im} \> \tau <1$ for $t \in (0,1)$. The analytic property of $f(z)$ is described by the following theorems :
\begin{theorem}[Theorem 1 in ~\cite{R}]
The function $f(z)$ is analytic on the upper half-plane $\mathbb{H}$ and for $z \in \mathbb{H}$ we have 
\begin{equation*}
   2\pi if(z) = f_1 (z) + f_2 (z) - e^{\frac{5}{2}z} \sum_{n=1}^\infty \frac{\varphi(n)}{n^\frac{5}{2}(z-\log n)}, \tag{1.4.2}
\end{equation*}where the last term on the right is meromorphic function on the whole complex plane with the poles at $z=\log n, \> n=1,2,\ldots$. The function
\begin{equation*}
   f_1 (z) = \int_{-\frac{1}{4}+i\infty}^{-\frac{1}{4}} \frac{\zeta(s-1)}{\zeta(s)}e^{sz}ds   \tag{1.4.3}
\end{equation*}is analytic on $\mathbb{H}$ and
\begin{equation*}
   f_2 (z) = \int_{\ell(-\frac{1}{4},\frac{5}{2})} \frac{\zeta(s-1)}{\zeta(s)}e^{sz}ds   \tag{1.4.4}
\end{equation*}is analytic on the whole complex plane.
\end{theorem}
\begin{theorem}[Theorem 2 in ~\cite{R}]
The function $f(z)$ can be continued analytically to a meromorphic function on the whole complex plane, which satisfies the functional equation
\begin{equation*}
   f(z) + \overline{f(\bar{z})} = B(z)   \tag{1.4.5}
\end{equation*}and
\begin{equation*}
   B(z) = -\frac{6}{\pi^2}e^{2z} + \frac{1}{2\pi^2}\sum_{k,n=1}^\infty \frac{\mu(k)}{n^2 k}\left[ \frac{1}{(nke^z -1)^2} + \frac{2}{nke^z -1} + \frac{1}{(nke^z + 1)^2} -\frac{2}{nke^z +1} \right],   \tag{1.4.6}
\end{equation*}where $B(z)$ is meromorphic function on the whole complex plane with the poles of the second order at $z=-\log nk, \> n,k=1,2,\ldots$. The only singularities of $f(z)$ are simple poles at the points $z=\log n \> (n=1,2,\ldots)$ on the real axis with residue 
\[ \res_{z=\log n} f(z) = -\frac{\varphi(n)}{2\pi i}, \] and the poles of the second order at $z=-\log m \> (m=1,2, \ldots)$ with residue
\[ \displaystyle \res_{z=-\log m} f(z) = \frac{1}{4\pi^2 m^2}\sum_{l|m}\mu(l)l. \] 
\end{theorem}

\subsection{The associated Euler totient function}

\hspace{1cm}

{\ka} defined the associated Euler totient function for a class of generalized $L$-functions including the Riemann zeta function, Dirichlet $L$-functions and obtained an asymptotic formula (see~\cite{Kac}) : By a polynomial Euler product we mean a function $F(s)$ of a complex variable $s=\sigma+it$ which for $\sigma>1$ is defined by the product of the form
\begin{equation*}
   F(s) = \prod_{p} F_p (s) = \prod_{p}\prod_{j=1}^d \left( 1-\frac{\alpha_j (p)}{ p^s } \right)^{-1},   \tag{1.5.1}
\end{equation*}where $p$ runs over primes and $|\alpha_j (p)| \leq 1$ for all $p$ and $1\leq j\leq d$. We assume that $d$ is chosen as small as possible, i.e. that there exists at least one prime number $p_0$ such that
\[ \ds \prod_{j=1}^d \alpha_j ( p_0 ) \neq 0. \] Then $d$ is called the \textit{Euler degree} of $F$. Note that the $L$-functions from number theory including the Riemann zeta function, Dirichlet $L$-functions, Dedekind zeta and Hecke $L$-functions of algebraic number fields, as well as the (normalized) $L$-functions of holomorphic modular forms and, conjecturally, many other $L$-functions are polynomial Euler products. For $F$ in (1.5.1) we define \textit{the associated Euler totient function} as follows :
\begin{equation}
   \varphi(n,F) = n\prod_{p|n} F_p (1)^{-1} \quad (n \in \mathbb{N}).    \tag{1.5.2}
\end{equation} Let
\begin{gather*}
   \gamma(p) = p\left( 1-\frac{1}{F_p (1)} \right),   \tag{1.5.3}
   \\
   C(F) = \frac{1}{2}\prod_p \left( 1 - \frac{\gamma(p)}{p^2} \right),   \tag{1.5.4}
\end{gather*}and
\begin{equation*}
   \alpha(n) = \mu(n)\prod_{p|n}\gamma(p).   \tag{1.5.5}
\end{equation*} By (1.5.2) and (1.5.3), we see that the Euler totient function $\varphi(n)$ and the twisted Euler $\varphi$-function $\varphi(n,\chi)$ correspond to the cases where $F$ is the Riemann zeta function $\zeta(s)$ and the Dirichlet $L$-function $L(s,\chi)$ respectively.
\begin{theorem}[Theorem 1.1 in ~\cite{Kac}]
For a polynomial Euler product $F$ of degree $d$ and $x \geq 1$ we have
\begin{equation*}
   \sum_{n \leq x} \varphi(n,F) = C(F)x^2 + O(x(\log 2x)^d).   \tag{1.5.6}
\end{equation*}
\end{theorem}
\begin{remark}Let us observe that $\alpha(n) \ll n^\epsilon$ for every positive $\epsilon$. Hence the series
\begin{equation}
   \sum_{n=1}^\infty \frac{\alpha(n)}{n^s}   \tag{1.5.7}
\end{equation}absolutely converges for $\sigma>1$ (see p33 in ~\cite{Kac}). Also, $\alpha(n)$ is multiplicative by (1.5.5). Therefore, 
\begin{equation*}
   \sum_{n=1}^\infty \frac{\alpha(n)}{n^2}
    = 2C(F).   \tag{1.5.8}
\end{equation*}
\end{remark}

\begin{lemma}[Lemma 2.2 in ~\cite{Kac}]
The series
\begin{equation*}
   \sum_{n=1}^\infty \frac{\varphi(n,F)}{n^s}   \tag{1.5.9}
\end{equation*}converges absolutely for $\sigma>2$ and in this half-plane we have
\begin{equation*}
   \sum_{n=1}^\infty \frac{\varphi(n,F)}{n^s} = \zeta(s-1)\sum_{n=1}^\infty \frac{\alpha(n)}{n^s}.  \tag{1.5.10}
\end{equation*}In particular,
\begin{equation*}
   \varphi(n,F) = n\sum_{m|n}\frac{\alpha(m)}{m}.   \tag{1.5.11}
\end{equation*}
\end{lemma}

\begin{lemma}[Lemma 2.3 in ~\cite{Kac}]
For $\sigma>1$ we have
\begin{equation*}
   \sum_{n=1}^\infty \frac{\alpha(n)}{n^s} = \frac{H(s)}{F(s)},   \tag{1.5.12}
\end{equation*}where 
\begin{equation*}
   H(s) = \sum_{n=1}^\infty \frac{h(n)}{n^s}   \tag{1.5.13}
\end{equation*}converges absolutely for $\sigma>1/2$. Moreover, as $n$ runs over square-free positive integers we have
\begin{equation}
   h(n) \ll \frac{1}{n}\exp\left( c\frac{\log n}{\log \log (n+2)} \right),   \tag{1.5.14}
\end{equation}where $c$ is a positive constant which may depend on $F$ and other parameters. In particular for such $n$, $h(n)$ is bounded.
\end{lemma}

\begin{lemma}[Lemma 2.4 in ~\cite{Kac}]
Let $\alpha(n)$ be defined by (1.5.5). For $x \geq 1$
\begin{equation*}
   \sum_{n \leq x} \frac{|\alpha(n)|}{n} \ll (\log(2x))^d .  \tag{1.5.15}
\end{equation*}
\end{lemma} 
Now we provide the definition of the \textit{Selberg class} ${\se}$ for our later purpose as follows : $F \in {\se}$ if
\begin{enumerate}
\item[(i)](\textit{ordinary Dirichlet series}) $\displaystyle F(s) = \sum_{n=1}^\infty a_F (n)n^{-s}$, absolutely convergent for $\sigma > 1$; 
\item[(ii)](\textit{analytic continuation}) there exists an integer $m\geq0$ such that $(s-1)^m \cdot F(s)$ is an entire function of finite order; 
\item[(iii)](\textit{functional equation}) $F(s)$ satisfies a functional equation of type $\Phi(s) = \omega\overline{\Phi(1-\overline{s})}$, where
                                                           \begin{equation*} 
                                                              \Phi(s) = Q^s \prod_{j=1}^r \Gamma(\lambda_j s + \mu_j)F(s) = \gamma(s)F(s),   \tag{1.5.15}
                                                           \end{equation*} say, with $r\geq0, Q>0, \lambda_j >0$, Re\hspace{0.001cm} $\mu_j \geq 0$ and $|\omega| = 1$;
\item[(iv)](\textit{Ramanujan conjecture}) for every $\epsilon>0, a_F (n) \ll n^\epsilon$.  
\item[(v)](\textit{Euler product}) $\displaystyle F(s) = \prod_{p} \exp \left( \sum_{\ell=0}^\infty  \frac{b_F ( p^\ell )}{p^{\ell s}} \right)$, where $b_F (n) = 0$ unless $n=p^m$ with $m\geq1$, and $b_F (n) \ll n^\vartheta$ for some 
                                  $\vartheta<\frac{1}{2}$.
\end{enumerate}

Note that we understand an empty product is equal to $1$.

The aim of the present paper is firstly to obtain the generalization of Theorem 1.1.1 for a remainder term in an asymptotic formula of an arithmetic function which satisfies some special conditions. The method using there is applied to the remainder term in the asymptotic formula of the associated Euler totient function. Also, we obtain the generalization of Theorem1.1.2 for the remainder term in the asymptotic formula of the associated Euler totient function. Secondly, we obtain results similar to Theorem 1.4.1 and Theorem 1.4.2 for the generalized $L$-functions which satisfy some conditions.

\section{Main results}
\subsection{The generalization of Theorem1.1.1}

\begin{theorem}[Theorem in ~\cite{I}]
Let $\{ a (n) \}$ be a complex-valued arithmetical function for which the series
                                  \begin{equation*}
                                     \sum_{n=1}^\infty \frac{a (n)}{n^2}   \tag{2.1.1}
                                  \end{equation*}is convergent with the sum $2\alpha$, where $\alpha$ is a complex number. Let $\{ b (n) \}$ be the arithmetical function defined by
                                  \begin{equation}
                                     b (n) = \sum_{d|n} a (d)\frac{n}{d}.   \tag{2.1.2}
                                  \end{equation}Assume for $x$ tending to infinity
                                  \begin{equation*}
                                     \sum_{n \leq x} b (n) = M(x) + \text{Er}(x),   \tag{2.1.3} 
                                  \end{equation*}where
                                  \begin{gather*}
                                     M(x) := \alpha x^2,   \tag{2.1.4}
                                     \\
                                     \text{Er}(x) := \sum_{n \leq x} b (n) -M(x).   \tag{2.1.5}
                                  \end{gather*}Now, we consider the following Volterra integral equation of second type
                                   \begin{equation*}
                                      F_1 (x) - \int_{0}^x F_1 (t)\frac{dt}{t} = \text{Er}(x) \quad (x \geq 0).  \tag{2.1.6}
                                   \end{equation*}Then, for every complex number $A$, the function
                                   \begin{equation*}
                                      F_1 (x) = (f_1 (x) + A)x  \quad (x \geq 0),   \tag{2.1.7} 
                                   \end{equation*}is a solution of the integral equation (2.1.6) and these exhaust all solutions of (2.1.6). Here,
                                   \begin{equation*} 
                                      f_1 (x) = -\sum_{n=1}^\infty \frac{a(n)}{n}\left\{ \frac{x}{n} \right\}   \tag{2.1.8}
                                   \end{equation*}for every $x \geq 0$.
\end{theorem}
In ~\cite{Kac and Wie}, the arithmetical functions $a(n) \text{ and } b(n)$ correspond to $\mu(n) \text{ and } \varphi(n)$ respectively, and all the hypothesis are satisfied. As for the error term $\text{Er}(x)$, we have a bound similar to $\text{Er}(x) = o(x^2)$ as $x$ tends to infinity in mind. As usual, if we say a function $F_1 $ is a solution of (2.1.6), then we implicitly assume that the integral in (2.1.6) exists in the sense that the limit
\begin{equation}
   \lim_{\epsilon \to 0+} \int_{\epsilon}^x |F_1 (t)|\frac{dt}{t}   \tag{2.1.9}
\end{equation}exists. We use the same convention throughout this paper. The formula (2.1.7) is a generalization of the result of Theorem 1.1.1. Also, the function $f_1 (x)$ is locally bounded. In fact, by the condition of theorem 
\begin{align*}
   f_1 (x)
   &= -\sum_{n=1}^\infty \frac{a (n)}{n}\left\{ \frac{x}{n} \right\}
   \\
   &= -\sum_{n=1}^\infty \frac{a (n)}{n} \left( \frac{x}{n} - \left[ \frac{x}{n} \right] \right)
   \\
   &= -x\sum_{n=1}^\infty \frac{a (n)}{n^2} + \left( \sum_{n \leq x}\frac{a (n)}{n}\left[ \frac{x}{n} \right] + \sum_{n > x}\frac{a (n)}{n} \left[ \frac{x}{n} \right] \right)
   \\
   &= -2\alpha x + \sum_{n \leq x}\frac{a (n)}{n}\left[ \frac{x}{n} \right].
\end{align*}We generalize Theorem 1.1.1 and 1.1.2 for the remainder term in the asymptotic formula for the associated Euler totient function. For a polynomial Euler product $F$ of degree $d$, let us put
\begin{equation*}
   E(x,F) = \sum_{n \leq x} \varphi(n,F) - C(F)x^2 ,   \tag{2.1.10}
\end{equation*}and
 \begin{align*}
   f(x,F)
   &= -\sum_{n=1}^\infty \frac{\alpha(n)}{n}\left\{ \frac{x}{n} \right\}.   \tag{2.1.11}
\end{align*}In the same way as the proof of Theorem 2.1.1, we have the following Corollary. The following two results are due to the author, but have not been written before. 
\begin{cor}
The Volterra integral equation of second type
\begin{equation*}
   F_1 (x,F) - \int_{0}^x F_1 (t,F)\frac{dt}{t} = E(x,F) \quad (x \geq 0)  \tag{2.1.12}
\end{equation*}has the following solution 
\begin{equation*}
   F_1 (x,F) = (f(x,F) + A)x  \quad (x \geq 0)   \tag{2.1.13}
\end{equation*}for every complex number $A$ and these exhaust all solutions of (2.1.12).
\end{cor}
For $x\geq0$ let us put 
\begin{equation*}
   g(x,F) = \sum_{n=1}^\infty \alpha(n)\left( \left\{ \frac{x}{n} \right\}^{2} + \left[ \frac{x}{n} \right] \right).  \tag{2.1.14}
\end{equation*}In the same way as in the proof of Theorem1.1.2, the error term (2.1.10) can be splitted as follows.
\begin{theorem}
For $x \geq 1$ we have
\begin{equation*}
   E(x,F) = xf(x,F) + \frac{1}{2}g(x,F).   \tag{2.1.15}
\end{equation*}
\end{theorem}We consider the \textit{arithmetic part} and the \textit{analytic part} for $E(x,F)$ as follows :
\begin{equation*}
   E^{\text{AR}}(x, F) = xf(x, F) \quad \text{and} \quad E^{\text{AN}}(x,F) = \frac{1}{2}g(x,F)   \tag{2.1.16}
\end{equation*}with $f(x,F)$ and $g(x,F)$ given by (2.1.11) and (2.1.14) respectively. Corollary 2.1.2 and Theorem 2.1.3 have not seen published, but can be proved in the same way as in ~\cite{I} and ~\cite{i}.

\subsection{Results similar to Theorem1.4.1 and 1.4.2}

\hspace{1cm}

If a function $F \in \mathcal{S}$ has a polynomial Euler product (1.5.1), the subclass of $\mathcal{S}$ of the functions with polynomial Euler product is denoted by $\mathcal{S}^{\text{poly}}$. Secondly, we obtain results similar to Theorem 1.4.1 and Theorem 1.4.2 for a function $F$ belonging to $\mathcal{S}^{\text{poly}}$. Let $\rho$ denote the non-trivial zeros of $F$ with positive imaginary part. We assume that the order of $\rho$ is simple. Moreover, let $T_n$ denote a sequence of real numbers which yields appropriate grouping of the zeros which will be given later. For Im $z > 0$ and $F \in\mathcal{S}^{\text{poly}}$, we consider a function defined by
\begin{equation*}
   f(z,F) = \lim_{n \to \infty} \sum_{\substack{\rho \\ 0< \text{Im} \> \rho < T_n}} \frac{e^{\rho z}\zeta(\rho-1)}{F^\prime (\rho)}.   \tag{2.2.1}
\end{equation*}If there are trivial zeros of $F$ on the imaginary axis in $\mathbb{H}$, we incorporate into the summation. We see that the series in (2.2.1) converges (see Section 4.2).
\begin{definition}[p 34 in ~\cite{Kac}]
For $\sigma>1$ and $F \in {\s}$, we define the function $\mu_F$ as follows :
\begin{equation*}
   \frac{1}{F(s)} = \sum_{n=1}^\infty \frac{\mu_F (n)}{n^s} = \prod_p \prod_{j=1}^d \left( 1 - \frac{\alpha_j (p)}{p^s} \right).   \tag{2.2.2}
\end{equation*}
\end{definition}
\begin{remark}[p34 in ~\cite{Kac}]
By (2.2.2), $|\mu_F (n)| \leq \tau_d (n)$, where $\tau_d (n)$ is the divisor function of order $d$, so that $\zeta^d (s) = \sum_{n=1}^\infty \tau_d (n) / n^s$ for $\sigma>1$. In particular $\tau_1 (n) = 1$ for all $n$.
\end{remark}
Using (2.2.2), for $\sigma>2$
\begin{align*}
   \frac{\zeta(s-1)}{F(s)}
   &= \left( \sum_{l=1}^\infty \frac{\mu_F (l)}{l^s} \right)\left( \sum_{m=1}^\infty \frac{1}{m^{s-1}} \right) 
   \\
   &= \sum_{n=1}^\infty \frac{g(n)}{n^s},   \tag{2.2.3}
\end{align*}where
\begin{equation*}
   g(n) = \sum_{d|n} \mu_F (d)\frac{n}{d}.   \tag{2.2.4}
\end{equation*}
\begin{theorem}[Theorem 4.1 in ~\cite{iw}]
The function $f(z,F)$ is analytic on the upper half-plane $\mathbb{H}$ and for $z \in \mathbb{H}$ we have 
\begin{equation*}
   2\pi if(z,F) = f_1 (z,F) + f_2 (z,F) - e^{bz}\sum_{n=1}^\infty \frac{g(n)}{n^b (z-\log n)},   \tag{2.2.5}
\end{equation*}where the last term on the right is a meromorphic function on the whole complex plane with the poles at $z=\log n, \> n=1,2,\ldots$. The function 
\begin{equation*}
   f_1 (z,F) = \int_{a+i\infty}^a \frac{\zeta(s-1)}{F(s)}e^{sz}ds   \tag{2.2.6}
\end{equation*}is analytic on $\mathbb{H}$ and
\begin{equation*}
   f_2 (z,F) = \int_L \frac{\zeta(s-1)}{F(s)}e^{sz}ds   \tag{2.2.7}
\end{equation*}is analytic on the whole complex plane. The definition of $a,b$ and the path of integration $L$ in (2.2.7) are mentioned later in Section 4.2.
\end{theorem}

Now we assume that $(r, \lambda_j ) = (1,1)$ in the functional equation (1.5.15). The complex number $\mu_1$ when $r=1$ in (1.5.15) is hereafter referred to as $\mu$. 

\begin{theorem}[Theorem 4.2 in ~\cite{iw}]
For $F$ belonging to  ${\s}$ whose $(r,\lambda_j) = (1,1)$ in (1.5.15) and $0 \leq \mu < 1$, the function $f(z,F)$ has a meromorphic continuation to $y>-\pi$.
\end{theorem}The $L$-functions associated with holomorphic cusp forms and Dedekind zeta functions of the imaginary quadratic fields are examples of $F$ considering in Theorem 2.2.2. Let
\begin{equation*}
   \mathbb{H}^{-} = \{ z \in \mathbb{C} : \text{Im}\hspace{0.1cm}  z < 0 \}.   \tag{2.2.8}
\end{equation*}We consider the function for $z \in \mathbb{H}^-$ and $F \in \mathcal{S}^{\text{poly}}$
\begin{equation*}
   f^{-} (z,F) = \lim_{n \to \infty} \sum_{\substack{\rho \\ -T_n < \text{Im} \> \rho < 0}} \frac{e^{\rho z}\zeta(\rho-1)}{F^\prime (\rho)}.   \tag{2.2.9}
\end{equation*} If there are trivial zeros of $F$ on the imaginary axis in $\mathbb{H}^-$, we incorporate into the summation. The convergence for the series on the right hand side in (2.2.9) is proved by the same way as in section 4.2. 
\begin{cor}[Corollary 4.3 in ~\cite{iw}]
For $F$ belonging to ${\s}$ which satisfies the same condition as in Theorem 2.2.2, the function $f^{-} (z,F)$ has a meromorphic continuation to $y<\pi$. 
\end{cor}
\begin{theorem}[Theorem 4.4 in ~\cite{iw}]
For $F$ belonging to ${\s}$ which satisfies the same condition as in Theorem 2.2.2, the function (2.2.1) can be continued analytically on the whole complex plane. In addition to the condition as in Theorem 2.2.2, we assume that the coefficient $a_F (n)$ in the Dirichlet series is real value for all $n$. Then, the function (2.2.1) satisfies the functional equation
\begin{equation*}
   f(z,F) + \overline{f(\overline{z},F)} = B(z,F),  \tag{2.2.10}
\end{equation*}where
\begin{equation*}
   B(z,F) = \frac{1}{2\pi i}(f_1 (z,F) + f_{1}^{-}(z,F)) - \frac{e^{2z}}{F(2)}.   \tag{2.2.11}
\end{equation*}for all $z \in \mathbb{C}$. The definition of $f_{1}^- (z,F)$ is mentioned later (4.6.4). 
\end{theorem} If a function $F \in {\se}$ satisfies the conditions (i)-(iii) on ${\se}$, we denote the this class by ${\se}^\#$ and call the \textit{extended Selberg class}. It is known that the Dirichlet coefficient $a_F (n)$ of $F \in {\se}^{\#}$ which satisfies some special conditions is real (see~\cite{kp}). We aim to prove the estimate corresponding to Theorem 1.1.5 for the remainder term in the asymptotic formula of the associated Euler totient function. To do this, we will use Theorem 2.2.1,2,4 and Corollary 2.2.3. 

\subsection{Motivation for Theorem 2.2.1,2,4 and Corollary 2.2.3}

\hspace{1cm}

We imitate the proof of Theorem1.1.5.\hspace{-0.01cm}Then, we need the following lemma which is the generalization of Lemma1.2.1. 
\begin{lemma}
For $\sigma>2$ and $F \in {\s}$, we have
\begin{equation*}
   \int_{1}^\infty E^{\text{AN}}(x,F)x^{-s-1}dx = \frac{C(F)}{s-2} + \frac{\zeta(s-1)}{s(1-s)}\frac{H(s)}{F(s)},   \tag{2.3.1}
\end{equation*}where the function $H(s)$ is the same as in (1.5.13).
\end{lemma}

\textbf{Proof of Lemma 2.3.1}:According to (2.1.15) and (2.1.16), we have $E^{\text{AN}}(x,F) = E(x,F) - xf(x,F)$. By inserting this into the integral in (2.3.1), we have
\begin{equation*}
   \int_{1}^\infty  E^{\text{AN}}(x,F)x^{-s-1}dx
   = \int_{1}^\infty E(x,F)x^{-s-1}dx - \int_{1}^\infty f(x,F)x^{-s}dx.  \tag{2.3.2}
\end{equation*} By the asymptotic formula (1.5.6) and $\sigma>2$, we have
\begin{equation*}
   \int_{1}^\infty E(x,F)x^{-s-1}dx
   = \int_{1}^\infty A(x,F)x^{-s-1}dx - \frac{C(F)}{s-2},   \tag{2.3.3}
\end{equation*}where
\[ \displaystyle A(x,F) = \sum_{n \leq x} \varphi(n,F). \]When we calculate the integral on the right hand side in (2.3.3), we use the following lemma :
\begin{lemma}[Theorem 4.2 in ~\cite{a}(Abel's identity)]
For any arithmetical function $a(n)$ let
\[ \ds A(x) = \sum_{n \leq x} a(n), \]where $A(x) = 0$ if $x<1$. Assume $f$ has a continuous derivative on the interval $[y,x]$, where $0<y<x$. Then we have
\begin{equation*}
   \sum_{y < n \leq x} a(n)f(n) = A(x)f(x) - A(y)f(y) - \int_{y}^x A(t)f^\prime (t)dt.
\end{equation*}
\end{lemma}Using Lemma 2.3.2 and (1.5.6) again, the integral on the right hand side in (2.3.3) is
\begin{equation*}
   \int_{1}^\infty A(x,F)x^{-s-1}dx 
   = \frac{1}{s}\sum_{n=1}^\infty \frac{\varphi(n,F)}{n^s}.  \tag{2.3.4}
\end{equation*}Since $\sigma>2$, we can apply Lemma1.5.2 and Lemma1.5.3. Hence,
\begin{align*}
   \sum_{n=1}^\infty \frac{\varphi(n,F)}{n^s}
   &= \zeta(s-1)\sum_{n=1}^\infty \frac{\alpha(n)}{n^s}
   \\
   &= \zeta(s-1)\frac{H(s)}{F(s)}.
\end{align*}Therefore, we have
\begin{equation*}
   \int_{1}^\infty E(x,F)x^{-s-1}dx = \frac{\zeta(s-1)}{s}\frac{H(s)}{F(s)} - \frac{C(F)}{s-2}.  \tag{2.3.5}
\end{equation*}On the other hand, by (2.1.11),
\begin{align*}
   -\int_{1}^\infty f(x,F)x^{-s}dx
   &= \int_{1}^\infty \left( \sum_{n=1}^\infty \frac{\alpha(n)}{n}\left\{ \frac{x}{n} \right\} \right)x^{-s}dx
   \\
   &= \int_{1}^\infty \left( x\sum_{n=1}^\infty \frac{\alpha(n)}{n^2} - \sum_{n=1}^\infty \frac{\alpha(n)}{n}\left[ \frac{x}{n} \right] \right)x^{-s}dx.
\end{align*}Since the series of the first term on the above second line is absolutely convergent, using (1.5.11), (1.5.8),
\begin{align*}
   \int_{1}^\infty \left( x\sum_{n=1}^\infty \frac{\alpha(n)}{n^2} - \sum_{n=1}^\infty \frac{\alpha(n)}{n}\left[ \frac{x}{n} \right] \right)x^{-s}dx
   &= \left( \sum_{n=1}^\infty \frac{\alpha(n)}{n^2} \right)\int_{1}^\infty x^{-s-1}dx - \int_{1}^\infty B(x,F)x^{-s}dx
   \\
   &= \frac{2C(F)}{s-2} - \int_{1}^\infty B(x,F)x^{-s}dx,
\end{align*}where
\[ \ds B(x,F) = \sum_{n \leq x} \frac{\varphi(n,F)}{n}. \]Using Lemma 2.3.2 again, we have
\begin{align*}
   \int_{1}^\infty B(x,F)x^{-s}dx
   &= \frac{1}{s-1}\sum_{n=1}^\infty \frac{\varphi(n,F)}{n^s}
   \\
   &= \frac{\zeta(s-1)}{s-1}\frac{H(s)}{F(s)}.
\end{align*}Therefore, we have
\begin{equation*}
   \int_{1}^\infty f(x,F)x^{-s}dx = -\frac{2C(F)}{s-2} + \frac{\zeta(s-1)}{s-1}\frac{H(s)}{F(s)}.   \tag{2.3.6}
\end{equation*}By inserting (2.3.5) and (2.3.6) into the right hand side in (2.3.2),  Lemma 2.3.1 follows.   $\square$

By the definition (2.1.14) and Lemma1.5.4, we see that $x^{-\sigma-1}E^{\text{AN}}(x,F) \in L(1,\infty)$. Also, by (2.1.15) and the definition (2.1.11), we see that $E^{\text{AN}}(x,F)$ is of bounded variation. Hence we can take the inverse Mellin transform in (2.3.1) by Fact1.3.1, and we obtain 
\begin{align*}
   E^{\text{AN}}(x,F) 
   &= \frac{1}{2\pi i}\int_{3-i\infty}^{3+i\infty}\left\{ \frac{C(F)}{s-2} + \frac{\zeta(s-1)}{s(1-s)}\frac{H(s)}{F(s)} \right\}x^s ds
   \\
   &= \frac{C(F)}{2\pi i}\int_{3-i\infty}^{3+i\infty}\frac{x^s}{s-2}ds + \frac{1}{2\pi i}\int_{3-i\infty}^{3+i\infty}  \frac{\zeta(s-1)}{s(1-s)}\frac{H(s)}{F(s)}x^s ds.  \tag{2.3.7}
\end{align*}When we calculate the first integral on the above second line in (2.3.7), we use the following lemma on contour integrals.
\begin{lemma}[Lemma 4 in ~\cite{a}]
If $c>0$, then if $a$ is any positive real number, we have
\begin{equation*}
   \frac{1}{2\pi i}\int_{c-i\infty}^{c+i\infty} \frac{a^z}{z}dz 
   = \begin{cases}
         1 \quad \text{if} \quad a>1,
         \\
         \frac{1}{2} \quad \text{if} \quad a=1,
         \\
         0 \quad \text{if} \quad 0<a<1.
      \end{cases}
\end{equation*}
\end{lemma}Using Lemma 2.3.3 and the residue theorem, we have 
\begin{equation*}
   E^{\text{AN}}(x,F) = \frac{1}{2\pi i}\int_\mathcal{L} \frac{\zeta(s-1)}{s(1-s)}\frac{H(s)}{F(s)}x^s ds,   \tag{2.3.8} 
\end{equation*}where the path of integration $\mathcal{L}$ is the same as in (1.3.1). We split the integral of left hand side in (2.3.8) into three parts
\begin{align*} 
   E^{\text{AN}}(x,F) 
   &= \frac{1}{2\pi i}\int_\mathcal{L} \zeta(s-1)\frac{H(s)}{F(s)} \left\{ -\frac{1}{s^2} - \frac{1}{s^3} + \frac{1}{s^3 (1-s)} \right\} x^s ds
   \\
   &= -\frac{1}{2\pi i} \int_\mathcal{L} \zeta(s-1)\frac{H(s)}{F(s)} \frac{x^s}{s^2}ds - \frac{1}{2\pi i} \int_\mathcal{L} \zeta(s-1)\frac{H(s)}{F(s)} \frac{x^s}{s^3}ds 
   \\
   &+  \frac{1}{2\pi i} \int_\mathcal{L} \zeta(s-1)\frac{H(s)}{F(s)} \frac{x^s}{s^3 (1-s)}ds.  \tag{2.3.9} 
\end{align*}Since the series (1.5.9) converges absolutely for $\sigma>2$, we can change the order of summation and integration. By calculating the residue at $s=0$, we have by (1.1.12)
\begin{equation*}
   \frac{1}{2\pi i} \int_\mathcal{L} \zeta(s-1)\frac{H(s)}{F(s)} \frac{x^s}{s^2}ds = R_1 (x,F)   \tag{2.3.10}
\end{equation*}and
\begin{equation*}
   \frac{1}{2\pi i} \int_\mathcal{L} \zeta(s-1)\frac{H(s)}{F(s)} \frac{x^s}{s^3}ds  = R_2 (x,F).  \tag{2.3.11}
\end{equation*}We use the following fact to obtain (2.3.10) and (2.3.11) :
\begin{fact}[Riesz typical means,~\cite{mv}]
For positive integers $k$ and positive real $x$ put
\begin{equation*}
   R_k (x) = \frac{1}{k!}\sum_{n \leq x} a_n \left( \log (x/n) \right)^k.  \tag{2.3.12}
\end{equation*}Then \[ R_k (x) = \int_{0}^x R_{k-1}(u)\frac{du}{u}, \]where
\[ \ds R_0 (x) = A(x) = \sum_{n \leq x} a_n . \]
\end{fact}When we calculate the last integral in (2.3.9), we have to consider the sum, corresponding to the sum in (1.3.4), that is
\begin{equation*}
   \sum_{\rho} \zeta(\rho-1)\frac{H(\rho)}{F^\prime(\rho)}\frac{x^\rho}{\rho^3 (1-\rho)}, \tag{2.3.12}
\end{equation*}where the summation is over non-trivial zeros of $F$. However, there are two problems on the sum (2.3.12). First, we do not know the behavior of $H(s)$ for $\sigma \leq \frac{1}{2}$. Therefore, We do not know whether $H(\rho)$ is defined and so, we do not assume the GRH i.e. $F(s) \neq 0 \> \left( \sigma> \frac{1}{2} \right)$. Secondly, we can not shift the line of integration $\mathcal{L}$ to the left from the critical line $\sigma=\frac{1}{2}$. If $F$ is the Riemann zeta function $\zeta(s)$ in (1.5.1), then by (1.5.2) the associated Euler totient function $\varphi(n,F)$ corresponds to the Euler totient function $\varphi(n)$. Since the Euler degree of $\zeta(s)$ equals $1$, we have $\gamma(p) = 1$ in (1.5.3). By (1.5.5), we have $\alpha(n) = \mu(n)$. By (1.5.12), $H(s) = 1$ for all $s$. When $H(s)=1$,  the sum (2.3.12) corresponds to the sum on the right hand side in (1.3.4). Therefore, the above problems are dissolved. To extend Theorem 1.1.5 for $F \in {\s}$, we have to consider the function for $z \in \mathbb{H}$
\begin{equation*}
   s(z, F) = \sum_{\rho} \frac{H(\rho)}{F^\prime (\rho)}\zeta(\rho-1)e^{\rho z}   \tag{2.3.13}
\end{equation*}which is a generalization for (1.4.1). However, the problem on the behavior of $H(s)$ for $\sigma \leq \frac{1}{2}$ is remained. To avoid these problems, we consider the sum of the case $H(s) = 1$, that is, the function (2.2.1).We need to generalize Theorem 1.4.1and Theorem 1.4.2 to prove the generalization of Theorem 1.1.5 for $F \in {\s}$ which satisfies some conditions. Obtaining the $\Omega$-estimate of $E^{\text{AN}}(x,F)$ for $F \in {\s}$ connects with the generalization of Theorem1.1.5. Conjecturally, every $F \in {\se}$ has an Euler product of type (1.5.1) and satisfies the GRH (see~\cite{Kac}).Therefore, proving the generalization of Theorem 1.1.5 is to obtain the $\Omega$-estimate of $E^{\text{AN}}(x,F)$ for many $L$-functions including the Riemann zeta function and is significant. That is why it is also significant to generalize Theorem 1.4.1 and Theorem 1.4.2 which will be necessary to prove the generalization of Theorem 1.1.5. In this paper, we could not obtain complete generalizations of Theorem 1.4.1 and Theorem 1.4.2. However, we could obtain results similar to Theorem 1.4.1 and Theorem 1.4.2 for $F \in {\s}$ which satisfy some special conditions. We want to prove the generalization of Theorem1.1.5 primarily. However, we try to obtain similar result to Theorem1.1.5 using these results. 

\section{Proof of Theorem 2.1.1}

\subsection{Preliminaries}

\hspace{1cm}

We prove Theorem 2.1.1. We define the auxiliary function for $x \geq 0$ by
\begin{equation*}
   R(x) = \text{Er}(x) - xf_1  (x).   \tag{3.1.1}
\end{equation*}First, we prepare the following two lemmas.
\begin{lemma}[LEMMA 1 in ~\cite{I}]
For all positive $x$,
                         \begin{equation*}
                            R(x) = -\int_{0}^x f_1 (t)dt,  \tag{3.1.2}
                         \end{equation*}where the function $f_1 (t)$ is the same as in (2.1.8).
\end{lemma} 
\textit{Proof.} Let us observe that $R(x)$ is a continuous function. For $x = 0$ and for positive $x$ which is not an integer, it is evident. Let $N$ be a positive integer. By splitting the series (2.1.8) at $N$, and considering the limit $\left\{ (N+x)/n \right\}$ as $x$ tending to $0$, we see that
\begin{align*}
   f_1 (N+0)
   &= -\sum_{n=1}^\infty \frac{a (n)}{n}\left\{ \frac{N+0}{n} \right\},
   \\
   f_1 (N-0)
   &= -\sum_{n=1}^\infty \frac{a (n)}{n}\left\{ \frac{N-0}{n} \right\}.
\end{align*}Since
\begin{equation*}
   \left\{ \frac{N+0}{n} \right\} - \left\{ \frac{N-0}{n} \right\}
   =\begin{cases}
        0  \quad \text{($n \nmid N$)},
        \\
        -1 \quad \text{($n \mid N$)}
    \end{cases}
\end{equation*}(see~\cite{Kac and Wie}, P2691), we have
\[ \ds f_1 (N+0) - f_1 (N-0) = \sum_{n|N} \frac{a (n)}{n} = \frac{b (N)}{N}. \]
Therefore
\begin{eqnarray*}
   R(N+0) - R(N-0)
\hspace{-0.25cm}   &=& \hspace{-0.25cm} (\text{Er}(N+0) - \text{Er}(N-0)) - N( f_1 (N+0) - f_1 (N-0) )
   \\
   &=& \hspace{-0.25cm} b (N) - N \cdot \frac{b (N)}{N}
   \\
   &=& \hspace{-0.25cm} 0,
\end{eqnarray*}and hence $R(N-0) = R(N+0) =R(N)$. Let $x$ be positive and not an integer. Take derivatives of the both sides of (3.1.1). Since $x$ is not a positive integer, we have $\text{Er}^\prime (x) = -M^\prime (x) = -2\alpha x$. Therefore we have
\[ R^\prime (x) = -2\alpha x - f_1 (x)-xf_{1}^\prime (x). \]For $x$ which is positive and not an integer, we have $\{ x/n \}^\prime = 1/n$ (see~\cite{Kac and Wie}, p2691). Considering the hypothesis on the series (2.1.1), differentiating term by term we obtain
\[ \ds f_{1}^\prime (x) = -\sum_{n=1}^\infty \frac{a (n)}{n} \cdot \frac{1}{n} = -2\alpha. \]Consequently, we have
\[ R^\prime (x) = -f_{1}(x) \]for $x$ which is positive and not an integer. Because of the fact $R(0) = 0$ and the continuity of $R(x)$, we have (3.1.2) for all positive $x$. $\square$
\begin{lemma}[LEMMA 2 in ~\cite{I}]
Let $G$ be a complex-valued function defined on $[0, \infty)$ satisfying
                         \begin{equation}
                            \int_{0}^x |G(t)|\frac{dt}{t} < +\infty   \tag{3.1.3}
                         \end{equation}and the integral equation
                         \begin{equation}
                            G(x) - \int_{0}^x G(t)\frac{dt}{t} = 0   \tag{3.1.4}
                         \end{equation}for all $x \geq 0$. Then we have
                         \begin{equation}
                            G(x) = Ax    \tag{3.1.5}
                        \end{equation}for some complex number $A$.
\end{lemma}
\textit{Proof.} It is obvious that (3.1.5) satisfies (3.1.3) and (3.1.4) for all $x \geq 0$. Conversely, take a function $G(x)$ arbitrarily satisfying (3.1.3) and (3.1.4) for all $x \geq 0$. By (3.1.3) and (3.1.4), we see that 
\[ G(x) = \int_{0}^x G(t)\frac{dt}{t} \] is a continuous function on $[0,+\infty)$. Thus, using integral equation again and using the fundamental theorem of calculus, we see that $G(x)$ is continuously differentiable on $(0,+\infty)$. By taking the derivative of (3.1.4), we have
\[ G^\prime (x) = \frac{G(x)}{x} \quad (x>0). \]Thus, we have $G(x) = Ax$ for $x > 0$ for some $A$ and by the continuity this holds for $x \geq 0$. $\square$

\subsection{Proof of Theorem 2.1.1} 

\hspace{1cm}

Let a function $F_1 (x)$ be the solution of the Volterra integral equation of second type (2.1.6) satisfying the condition (2.1.9). Using (2.1.6) and (3.1.2), from (3.1.1) we have 
\begin{equation*} 
   \int_{0}^x t^{-1}( F_1 (t) - tf_1 (t) )dt = F_1 (x) - xf_1 (x) \quad (x \geq 0).   \tag{3.2.1}
\end{equation*}Now we put 
\begin{equation*}
   G(x) := F_1 (x) - xf_1 (x).   \tag{3.2.2} 
\end{equation*}Then, the equation (3.2.1) yields
\begin{equation*} 
   \int_{0}^x t^{-1}G(t)dt = G(x) \quad (x \geq 0).    \tag{3.2.3}
\end{equation*}Using Lemma 3.1.2, we must have (3.1.5). By substituting into (3.2.2), we have the solution (2.1.7). Conversely, if we assume that $F_1 (x)$ is a function of type (2.1.7). Then,
\begin{align*}
   F_1 (x) - \int_{0}^x F_1 (t)\frac{dt}{t}
   &= (f_1 (x) + A)x - \int_{0}^x (f_1 (t) + A)dt
   \\
   &= (f_1 (x) + A)x - \int_{0}^x f_1 (t)dt - Ax
   \\
   &= xf_1 (x) - \int_{0}^x f_1 (t)dt.
\end{align*}Using (3.1.2) and (3.1.1), 
\begin{align*}
   xf_1 (x) - \int_{0}^x f_1 (t)dt
   &= xf_1 (x) + R(x)
   \\ 
   &= xf_1 (x) +  \text{Er}(x) - xf_1 (x)
   \\
   &= \text{Er}(x).
\end{align*}Therefore, the function $F_1 (x)$ of type (2.1.7) is the solution of the integral equation (2.1.6) for all $x \geq 0$. Since the function $f_1 (x)$ is a locally bounded as noted in section 2.1, and $A$ is a constant, it is clear that the function $F_1 (x)$ satisfies the condition (2.1.9). The completes the proof. $\square$

\section{Proof of Theorem 2.2.1,2,4 and Corollary 2.2.3}
\subsection{Some auxiliary results on the Whittaker function}

\hspace{1cm}

First, we recall the definition of the Whittaker function which is necessary to show the main theorems.To do this,we introduce some related functions. Secondly, we prepare some auxiliary results e.g.\hspace{-0.0005cm}the integral expression and the asymptotic expansion.
\begin{definition}[The confluent hypergeometric function (~\cite{L})]
Let $z$ be a complex variable, $\alpha$ and $\gamma$ are parameters which can take arbitrary real or complex values except that $\gamma \neq 0, -1, -2, \ldots$. Moreover,
\begin{equation*}
   (\lambda)_0 = 1, \quad (\lambda)_k = \frac{\Gamma(\lambda+k)}{\Gamma(\lambda)} = \lambda(\lambda+1)\cdots (\lambda+k-1) \quad (k \in \mathbb{N}).  \tag{4.1.1}
\end{equation*}We define the \textit{confluent hypergeometric function (Kummer's function)} as follows :
\begin{equation*} 
\Phi(\alpha, \gamma ; z) = \sum_{n=0}^\infty \frac{(\alpha)_n}{(\gamma)_n}\frac{z^n}{n!} \quad (|z|<\infty).  \tag{4.1.2}
\end{equation*} 
\end{definition}By ratio test, the series (4.1.2) is convergence absolutely for all $\alpha,\gamma$ and $z$, except $\gamma \neq 0,-1,-2,\ldots$. Hence, (4.1.2) is an analytic and one-valued function for all $z$. It is also to be noted that (4.1.2) is a particular solution of the linear differential equation ( \textit{Kummer's equation} )
\begin{equation*}
   z\frac{d^2 u}{dz^2} + (\gamma-z)\frac{du}{dz} - \alpha u = 0,   \tag{4.1.3}
\end{equation*}where $ \alpha, \gamma$ are the same as in (4.1.2).
\begin{definition}[The confluent hypergeometric function of the second kind (~\cite{L})]
We introduce a new function
\begin{gather*}
   \Psi(\alpha, \gamma ; z) = \frac{\Gamma(1-\gamma)}{\Gamma(1+\alpha-\gamma)}\Phi(\alpha,\gamma;z) + \frac{\Gamma(\gamma-1)}{\Gamma(\alpha)}z^{1-\gamma}\Phi(1+\alpha-\gamma,2-\gamma;z),
   \\
  (|\arg z| < \pi, \gamma \neq 0, \pm 1, \pm 2, \ldots)  \tag{4.1.4}
\end{gather*}called \textit{the confluent hypergeometric function of the second kind}. The condition $\gamma \neq 0, \pm 1, \pm 2, \ldots$ in (4.1.4) comes from the condition of the $\Gamma$-factors in the numerator, $\Phi(\alpha, \gamma ; z)$ and $\Phi(1+\alpha-\gamma,2-\gamma;z) $ on the right hand side in (4.1.4).
\end{definition}Since the function (4.1.4) is a many-valued function of $z$ for $\alpha$ and $\gamma$ real or complex, we take its branch which lies in the $z$-plane cut along the negative real axis. Also, (4.1.4) is analytic function for all $\alpha, \gamma$ and $z$ except $\gamma \neq 0,-1,-2,\ldots$.
\begin{definition}[The Whittaker function (~\cite{L})]
A class of functions related to the confluent hypergeometric functions, and often encountered in the applications, consists of the \textit{Whittaker function}, defined by the formula  
\begin{equation*}
   W_{k, \hspace{0.01cm} \mu}(z) = z^{\mu+\frac{1}{2}}e^{-\frac{z}{2}}\Psi\left( \frac{1}{2}-k+\mu, 2\mu+1 ; z \right) \quad (|\arg z| < \pi).  \tag{4.1.5}
\end{equation*}
\end{definition}By (4.1.4) or (4.1.5), it follows that (4.1.5) is a many-valued function of $z$. Therefore, we also take its branch as the same in (4.1.4).
\begin{theorem}[Barnes type integral for the Whittaker function (~\cite{Sl}, ~\cite{WW})]
The Barnes integral for $W_{k,\mu} (z)$ is
\begin{equation*}
W_{k,\mu}(z) = \frac{e^{-\frac{1}{2}z} z^k}{2\pi i}\int_{c-\infty i}^{c+\infty i}\frac{\Gamma(s)\Gamma\left( -s-k-\mu+\frac{1}{2} \right)\Gamma\left( -s-k+\mu+\frac{1}{2} \right)}{\Gamma\left( -k-\mu+\frac{1}{2} \right)\Gamma\left( -k+\mu+\frac{1}{2} \right)}z^s ds,   \tag{4.1.6}
\end{equation*}where $|\arg z| < \frac{3}{2}\pi$, and neither of the numbers $k\pm \mu + \frac{1}{2}$ is a positive integer or zero ; the contour has deformations if necessary so that the poles of $\Gamma(s)$ and those of $\Gamma\left( -s-k-\mu +\frac{1}{2} \right) \Gamma\left( -s-k+\mu +\frac{1}{2} \right)$ are on opposite sides of it. 
\end{theorem}
In (4.1.6), it holds for all finite values of $c$ provided that the contour of integration can always be deformed so as to separate the poles $\Gamma(s)$ and those of the other $\Gamma$-factors.The integral in (4.1.6) represents a function of $z$ which is analytic at all points in the domain $|\arg z| \leq \frac{3}{2}\pi - \alpha$, where $\alpha$ is any positive number using the following fact.

\begin{fact}[Stirling's formula,~\cite{WW}]
For $s=\sigma + it$, 
\begin{equation*}
   \Gamma(s) \ll e^{-\frac{\pi}{2}|t|}|t|^{\sigma-\frac{1}{2}}
\end{equation*}as $|t| \to \infty$.
\end{fact}
\begin{theorem}[The asymptotic expansions in $z$ for $\Psi(a;b;z)$ (~\cite{Sl})]
We find that, as $z \to 0$,
\begin{align*}
   \Psi(a;b;z)
   &= \frac{\Gamma(b-1)}{\Gamma(a)}z^{1-b}+O(|z|^{\text{Re} b -2}) \quad (\text{Re} \hspace{0.1cm} b \geq 2, b \neq 2),   \tag{4.1.7}
   \\
   &= \frac{\Gamma(b-1)}{\Gamma(a)}z^{1-b}+O(|\log z|) \quad (b = 2),   \tag{4.1.8}
  \\
   &= \frac{\Gamma(b-1)}{\Gamma(a)}z^{1-b}+O(1) \quad (1 < \text{Re} \hspace{0.1cm} b < 2),   \tag{4.1.9}
   \\
   &= \frac{\Gamma(1-b)}{\Gamma(1+a-b)} + \frac{\Gamma(b-1)}{\Gamma(a)}z^{1-b} + O(|z|) \quad (\text{Re} \hspace{0.1cm} b =1, b \neq 1),   \tag{4.1.10}
   \\
   &= -\frac{1}{\Gamma(a)}\left\{ \log z + \frac{\Gamma^\prime}{\Gamma}(a) + 2C_0 \right\} + O(|z\log z|) \quad (b=1),   \tag{4.1.11} 
\end{align*}where $C_0$ is Euler's constant.
\end{theorem}
By the definition (4.1.5) and Theorem 4.1.2, we have the following asymptotic expansions in $z \to 0$ for $W_{k,\mu}(z)$. 
\begin{theorem}[The asymptotic expansions in $z \to 0$ for $W_{k,\mu} (z)$]
\begin{align*}
   W_{k,\mu}(z) 
   &= \frac{\Gamma(2\mu)}{\Gamma\left( \frac{1}{2} + \mu - k \right)}z^{\frac{1}{2}-\mu}+O(z^{\frac{3}{2}- \text{Re} \mu}) \quad \left( \text{Re} \hspace{0.1cm} \mu \geq \frac{1}{2}, \mu \neq \frac{1}{2} \right),  \tag{4.1.12}
   \\
   &= \frac{1}{\Gamma(1-k)} + O(|z\log z|) \quad \left( \mu = \frac{1}{2} \right),  \tag{4.1.13}
   \\
   &= \frac{\Gamma(2\mu)}{\Gamma\left( \frac{1}{2} + \mu - k \right)}z^{ \frac{1}{2} - \mu} + O(|z|^{\text{Re}\mu + \frac{1}{2}})  \quad \left( 0 < \text{Re} \hspace{0.1cm} \mu < \frac{1}{2} \right),  \tag{4.1.14}
   \\
   &= \frac{\Gamma(-2\mu)}{\Gamma\left( \frac{1}{2} - \mu -k \right)}z^{\mu+\frac{1}{2}} + \frac{\Gamma(2\mu)}{\Gamma\left( \mu+\frac{1}{2} - k \right)}z^{-\mu+\frac{1}{2}} + O(|z|^{\text{Re}\mu + \frac{3}{2}}) \quad 
   ( \text{Re} \hspace{0.1cm} \mu = 0, \mu \neq 0 ),  \tag{4.1.15}
   \\
   &= -\frac{z^\frac{1}{2}}{\Gamma\left( \frac{1}{2} - k \right)}\left( \log z + \frac{\Gamma^\prime}{\Gamma}\left( \frac{1}{2} - k \right) + 2C_0 \right) + O(|z|^\frac{3}{2} |\log z|) \quad (\mu = 0).   \tag{4.1.16}
\end{align*}
\end{theorem}

\subsection{The convergence of $f(z,F)$}

\hspace{1cm}

In this section, we see that the series in (2.2.1) converges.We use the following lemma.
\begin{lemma}[Lemma 4. in ~\cite{Sr}]
Let $F \in \mathcal{S}$ and let $T$ be sufficiently large. Moreover, let $H=D\log\log T$ be fixed, where $D$ is a large positive constant. In any subinterval of length $1$ in $[T-H, T+H]$ there are lines $t=t_0$ such that
\begin{equation*}
   |F(\sigma+it_0)|^{-1} = O(\exp(C(\log T)^2)),   \tag{4.2.1}
\end{equation*}uniformly in $\sigma\geq -2$, where $C$ is a positive constant.
\end{lemma}Let $T$ be sufficiently large.We fix $H=D\log\log T$, where $D$ is a large positive constant.We take any subinterval $[n, n+1]$, where $n$  is a positive integer in $[T-H, T+H]$. Then, by Lemma 4.2.1 there are lines $t=T_n$ such that
\begin{align*}
|F(\sigma+iT_n)|^{-1} = O( \exp(C_1 (\log T)^2) )   \tag{4.2.2}
\end{align*}uniformly for $\sigma \geq -2$, where $C_1$ is a positive constant. Since $T_n$ is contained in the interval $[T-H,T+H]$, we can see that $T_n  \sim T$ as $n$ tends to infinity. Let $\alpha = \frac{1}{2}\min\{\text{Im}\hspace{0.1cm}\rho ; \text{Im}\hspace{0.1cm}\rho >0  \}$ and $\mathscr{L}$ denote the contour consisting of line segments
\[ \left[ b, b+i T_n \right], \> \left[ b+i T_n , a+i T_n \right], \> \left[ a+i T_n , a \right], \> \left[a,\frac{a+b}{2}+i\alpha \right], \> \left[\frac{a+b}{2}+i\alpha, b \right], \]where $\max \left\{ -\frac{3}{2}, \frac{1}{2}\max\left\{\text{Re} \hspace{0.1cm} \rho; \text{Re} \hspace{0.1cm} \rho < 0 \right\} \right\} < a < 0, b > \frac{5}{2}$. We assume that the real part of $s = a + it \hspace{0.1cm} (t \in \mathbb{R})$ does not coincide the poles of $\Gamma( s + \mu )\Gamma(s - \mu)$, where $0 \leq \mu < 1$. We consider the following contour integral round $\mathscr{L}$ :
\begin{equation*}
   \int_\mathscr{L} \frac{\zeta(s-1)}{F(s)}e^{sz}ds.  \tag{4.2.3}
\end{equation*} Since we assume the order of $\rho$ is simple, we have by residue theorem 
\begin{align*}
   \int_\mathscr{L} \frac{\zeta(s-1)}{F(s)}e^{sz}ds
   &= \int_{a+iT_n}^a \frac{\zeta(s-1)}{F(s)}e^{zs}ds + \int_L \frac{\zeta(s-1)}{F(s)}e^{zs}ds 
   \\
   &+ \int_b^{b+iT_n}\frac{\zeta(s-1)}{F(s)}e^{zs}ds + \int_{b+iT_n}^{a+iT_n}\frac{\zeta(s-1)}{F(s)}e^{zs}ds
   \\
   &= 2\pi i \sum_{\substack{\rho \\ 0< \text{Im} \> \rho < T_n}} \frac{e^{\rho z}\zeta(\rho-1)}{F^\prime (\rho)},   \tag{4.2.4}
\end{align*}where the path of integration $L$ consists of two line segments $\left[a,\frac{a+b}{2}+i\alpha \right]$ and $\left[\frac{a+b}{2}+i\alpha, b \right]$. There is a possibility that the zeros of $F$ lie on the interval $(0,1)$. That is why we consider the path of integration $L$. We estimate the integral along the line segment $[b+iT_n , a+iT_n]$. By (4.2.1), for $a \leq \sigma \leq b$, we have the estimate
\[ |F(\sigma+iT_n)|^{-1} = O( \exp(C(\log T)^2) ). \] For $z=x+iy \> (y>0)$,
\begin{align*}
   &\left| \int_{a+iT_n}^{b+iT_n} \frac{\zeta(s-1)}{F(s)}e^{zs}ds \right|
   \\
   &\leq \int_{a+iT_n}^{b+iT_n} \left| \frac{\zeta(s-1)}{F(s)}e^{zs} \right| |ds|
   \\
   &\ll \int_{a}^b |\zeta(\sigma-1+iT_n)|\exp( C(\log T)^2 + x\sigma -yT_n )d\sigma
   \\
   &\ll (b-a)\exp\{ C(\log T)^2 -yT_n + |x|(|a| + |b|)\} T_{n}^c ,   \tag{4.2.5}
\end{align*}where the constant $c$ may depend on $a,b$. The last term on the right hand side in the above tends to zero as $n$ tends to infinity. By Theorem 2.2.1, the convergence of the other integrals in (4.2.4) are ensured (see (2.2.5)-(2.2.7)). Therefore, the series in (2.2.1) is convergent. $\square$ 

\subsection{Proof of Theorem 2.2.1}

\hspace{1cm}

By (4.2.4), for $z \in \mathbb{H}$ we have
\[ f_1 (z,F) + f_2 (z,F) + f_3 (z,F) = 2\pi if(z,F), \] where $f_1 (z,F),  f_2 (z,F),  f_3 (z,F)$ denote corresponding integrals in (4.2.4), respectively. First, we calculate the integral
\[ f_3 (z,F) = \int_{b}^{b+i\infty} \frac{\zeta(s-1)}{F(s)}e^{zs}ds. \] Using the Dirichlet series expansion (2.2.2), we have
\begin{align*}
   f_3 (z,F)
   &= \int_{b}^{b+i\infty} \left( \sum_{n=1}^\infty \frac{g(n)}{n^s} \right)e^{zs}ds
   \\
   &= \sum_{n=1}^\infty g(n) \int_{b}^{b+i\infty}e^{s(z-\log n)}ds
   \\
   &= -e^{bz}\sum_{n=1}^\infty \frac{g(n)}{n^b (z-\log n)}.   \tag{4.3.1}
\end{align*}The interchange of the order of integration and summation is justified for $z \in \mathbb{H}$ by the absolute and uniform convergence of the series on the third line in (4.3.1). 

Secondly, we prove the function (2.2.7) is analytic on the whole complex plane.Since the length of $L$ is finite, $\zeta(s-1) = O(1)$. Also, there are no zeros on $L$, the function $\{F(s)\}^{-1}$ is bounded. Therefore, for $z=x+iy$ and $s=\sigma+it \> (a \leq \sigma \leq b, 0 \leq t \leq \alpha)$,
\begin{align*}
   \int_L \left| \frac{\zeta(s-1)}{F(s)}e^{zs} \right| |ds|
   &\ll \int_L \exp(x\sigma-yt)|\zeta(s-1)||ds|
   \\
   &= O\left( \int_L \exp(x\sigma-yt)|ds| \right).
\end{align*}In the case $y\geq 0$,
    \[ \int_L \exp(x\sigma-yt)|ds| \ll \int_L \exp(x\sigma)|ds| \ll_x 1. \]In the case $y<0$, 
    \[ \int_L \exp(x\sigma-yt)|ds| \ll_{x,y} 1. \]Hence the function (2.2.7) is analytic on the whole complex plane. 
    
Finally, we prove the function (2.2.6) is analytic on $\mathbb{H}$. Here, we use the following two lemmas related to the Selberg class.  
\begin{definition}[p29 in ~\cite{P}]
Let $F \in \mathcal{S}$, 
\begin{equation*}d_F = 2\sum_{j=1}^r \lambda_j   \tag{4.3.2}
\end{equation*} is the \textit{degree} of $F(s)$. 
\end{definition}
\begin{lemma}[Theorem 3.1 in ~\cite{C and G}]
If $F \in \mathcal{S}$, then $F=1$ or $d_F \geq 1$.
\end{lemma}
\begin{lemma}[~\cite{M}.(8), p423]
For $\sigma<0$,
\begin{equation*}
  |F(\sigma+it)| \asymp t^{d_F (\frac{1}{2}-\sigma)}|F(1-\sigma+it)|   \tag{4.3.3}
\end{equation*}as $t \to \infty$.
\end{lemma}Here, we use the notation$f \asymp g$, if both $f \ll g$ and $f \gg g$ hold. Using Lemma 4.3.2, we have
\begin{equation*}
     t^{\left( \frac{1}{2}-a \right)d_F} |F(1-a+it)| \ll |F(a+it)| \ll t^{\left( \frac{1}{2}-a \right)d_F}|F(1-a+it)|.
\end{equation*}Since $\max \left\{ -\frac{3}{2}, \frac{1}{2}\max\left\{\text{Re} \hspace{0.1cm} \rho; \text{Re} \hspace{0.1cm} \rho < 0 \right\} \right\} < a < 0$, so $\{F(1-a+it)\}^{-1}$ is bounded. Hence, we have
\[ |F(a+it)|^{-1} \ll t^{-(\frac{1}{2}-a)d_F}. \]
Also, by the functional equation for $\zeta(s)$,
\begin{align*}
   \zeta(s-1) 
   &= -\frac{1}{\pi}(2\pi)^{s-1}\cos\left( \frac{\pi}{2}s \right)\Gamma(2-s)\zeta(2-s)
   \\
   &\ll |t|^{\frac{3}{2}-a}.
\end{align*}In the case $F=1$, we have
\begin{align*}
   \int_{a+i\infty}^a \left| \frac{\zeta(s-1)}{F(s)}e^{zs} \right||ds|
   &\ll \int_{0}^\infty t^{\frac{3}{2}-a}e^{ax-yt}dt
   \\
   &\ll_{x,y} 1
   \\
   &\ll 1.
\end{align*}Next consider the case $F \neq 1$. Now, $F(s)$ belongs to $\mathcal{S}^{\text{poly}}$, so $F(s)$ is analytic except the point $s=1$. Since $\max \left\{ -\frac{3}{2}, \frac{1}{2}\max\left\{\text{Re} \hspace{0.1cm} \rho; \text{Re} \hspace{0.1cm} \rho < 0 \right\} \right\} < a < 0$, $|F(a+it)|^{-1}$ is bounded near $t=0$. Also, $d_F \geq 1$ by Lemma 4.3.1, we have 
\begin{align*}
   \int_{a+i\infty}^a \left| \frac{\zeta(s-1)}{F(s)}e^{zs} \right||ds|
   &\ll \int_{0}^\infty \frac{e^{ax-yt}}{|F(a+it)|}t^{\frac{3}{2}-a}dt
   \\
   &\ll \int_{0}^\infty t^{-\left( \frac{1}{2} - a \right)d_F + \frac{3}{2} - a }e^{ax-yt}dt
   \\
   &< \int_{0}^\infty t^{-\frac{1}{2}d_F + \frac{3}{2} - a }e^{ax-yt}dt.
\end{align*}Therefore, the integral on the above third line is absolutely and uniformly convergent on every compact subset of $\mathbb{H}$. Consequently, the function (2.2.6) is analytic for $y =$ Im $z >0$. $\square$

\subsection{Proof of Theorem 2.2.2}
We prove that the function $f(z,F)\>(z=x+iy)$ has a meromorphic continuation to $y>-\pi$. By Theorem 2.2.1, the function
\begin{align*} 
   f_1 (z,F) 
   &= \int_{a+i\infty}^a \frac{\zeta(s-1)}{F(s)}e^{zs}ds
   \\
   &= -\int_a^{a + i\infty} \frac{\zeta(s-1)}{F(s)}e^{zs}ds
\end{align*}is convergent for $y>0$. We recall the hypotheses that $(r, \lambda_j ) = (1,1)$ in (1.5.15) and $0 \leq \mu <1$. We rewrite the functional equation (1.5.15) under these hypotheses as follows :
\begin{align*}
   Q^s \Gamma( s + \mu )F(s) 
   &= \omega \overline{Q^{1-\overline{s}}\Gamma( 1 - \overline{s} + \mu )F(1-\overline{s})}
   \\
   &= \omega Q^{1-s}\Gamma(1 - s + \mu){\f},
\end{align*}where the conditions of $Q$ and $\omega$ are the same as noted in (1.5.15). Hence
\begin{equation*} 
   \frac{1}{F(s)} = \overline{\omega}Q^{2s-1}\frac{\Gamma( s + \mu )}{\Gamma( 1- s + \mu )}\frac{1}{{\f}}.   \tag{4.4.1}
\end{equation*}Using the following elementary formula for the $\Gamma$-function 
\[ \Gamma(s)\Gamma(1-s) = \frac{\pi}{\sin \pi s}, \](4.4.1) yields 
\begin{equation*}
   \frac{1}{F(s)} = \frac{\overline{\omega}}{\pi}Q^{2s-1}\sin\pi(s-\mu)\Gamma(s+\mu)\Gamma(s-\mu)\frac{1}{{\f}}.  \tag{4.4.2}
\end{equation*}By (4.4.2) and the functional equation for $\zeta(s)$, we have
\begin{align*}
   &f_1 (z,F)
   \\ 
   &=  -\int_a^{a + i\infty} \frac{\zeta(s-1)}{F(s)}e^{zs}ds
   \\
   &= \frac{\overline{\omega}}{2\pi^3 Q}\int_a^{a+i\infty} \hspace{-0.7cm} (2\pi Q^2)^s \cos\left( \frac{s}{2}\pi \right) \Gamma(2-s)\zeta(2-s)\sin\pi(s-\mu)
   \\
   &\times \Gamma(s+\mu)\Gamma(s-\mu)\frac{e^{zs}}{{\f}}ds
   \\
   &= \frac{\overline{\omega} e^{-\mu\pi i}}{(2\pi)^3 Qi}\int_a^{a+i\infty} (2\pi Q^2 )^s \frac{\zeta(2-s)}{{\f}}\Gamma(s - \mu)\Gamma(s + \mu )\Gamma(2-s)e^{\left( z + \frac{3}{2}\pi i \right)s}ds
   \\
   &-\frac{\overline{\omega} e^{\mu\pi i}}{(2\pi)^3 Qi}\int_a^{a+i\infty} (2\pi Q^2 )^s \frac{\zeta(2-s)}{{\f}}\Gamma( s - \mu )\Gamma(s + \mu)\Gamma(2-s)e^{\left(z - \frac{\pi}{2}i \right)s}ds   
   \\
   &+\frac{\overline{\omega} e^{-\mu\pi i} }{(2\pi)^3 Qi}\int_a^{a+i\infty} (2\pi Q^2 )^s \frac{\zeta(2-s)}{{\f}}\Gamma( s - \mu )\Gamma(s + \mu)\Gamma(2-s)e^{\left(z + \frac{\pi}{2}i \right)s}ds   
   \\   
   &-\frac{\overline{\omega} e^{\mu\pi i} }{(2\pi)^3 Qi} \int_a^{a+i\infty} (2\pi Q^2 )^s \frac{\zeta(2-s)}{{\f}}\Gamma(s - \mu )\Gamma( s + \mu )\Gamma(2-s) e^{(z - \frac{3}{2}\pi i)s}ds   
   \\
   &= f_{11}(z,F) + f_{12}(z,F) + f_{13}(z,F) + f_{14}(z,F),  \tag{4.4.3}
\end{align*}where $f_{11}(z,F), f_{12}(z,F), f_{13}(z,F), f_{14}(z,F)$ denote the corresponding integrals in (4.4.3) respectively. By Fact 4.1.2,
\begin{equation*}
   \Gamma(s+\mu)\Gamma(s-\mu)\Gamma(2-s) \asymp e^{-\frac{3}{2}\pi |t|}|t|^{a+\frac{1}{2}}
\end{equation*}as $t$ tends to infinity and $\max \left\{ -\frac{3}{2}, \frac{1}{2}\max\left\{\text{Re} \hspace{0.1cm} \rho; \text{Re} \hspace{0.1cm} \rho < 0 \right\} \right\} < a < 0$, $f_{11}(z,F)$ is analytic for $y > -3\pi$, $f_{12}(z,F)$ for $y> -\pi$, $f_{13}(z,F)$ for $y > -2\pi$, $f_{14}(z,F)$ for $y>0$. Splitting the integral in $f_{14}(z,F)$, we have
\begin{align*}
   f_{14}(z,F)
   &= -\frac{\overline{\omega} e^{\mu\pi i} }{(2\pi)^3 Qi} \int_a^{a+i\infty} \hspace{-0.5cm} (2\pi Q^2 )^s \frac{\zeta(2-s)}{{\f}}\Gamma(s - \mu )\Gamma( s + \mu )\Gamma(2-s) e^{(z - \frac{3}{2}\pi i)s}ds   
   \\   
   &= \frac{\overline{\omega}e^{\mu\pi i} i}{(2\pi)^3 Q}\left\{ \int_{a-i\infty}^{a+i\infty} - \int_{a-i\infty}^a \right\} (2\pi Q^2 )^s \frac{\zeta(2-s)}{{\f}}
   \\
   &\times \Gamma(s - \mu )\Gamma( s + \mu )\Gamma(2-s) e^{(z - \frac{3}{2}\pi i)s}ds.
\end{align*}We consider
\begin{equation*}
   I_1 (z,F) = \int_{a-i\infty}^{a+i\infty} (2\pi Q^2 )^s \frac{\zeta(2-s)}{{\f}}\Gamma(s - \mu )\Gamma( s + \mu )\Gamma(2-s) e^{(z - \frac{3}{2}\pi i)s}ds   \tag{4.4.4}
\end{equation*}and
\begin{equation*}
  I_2 (z,F) =  \int_{a-i\infty}^a (2\pi Q^2 )^s \frac{\zeta(2-s)}{{\f}}\Gamma(s - \mu )\Gamma( s + \mu )\Gamma(2-s) e^{(z - \frac{3}{2}\pi i)s}ds.   \tag{4.4.5}
\end{equation*}Here, we recall the hypothesis that the real part of $s=a+it \> (t \in \mathbb{R})$ does not coincide with the poles of $\Gamma(s+\mu)\Gamma(s-\mu)$, where $0 \leq \mu < 1$. By Fact 4.1.2, we can see that the integral $I_2 (z,F)$ is convergent for $y<3\pi$. Since the function $f_{14}(z,F)$ is analytic for $y>0$, the integral $I_1 (z,F)$ is convergent for $0<y<3\pi$. Using the Dirichlet series expansion (2.2.2), we have formally
\begin{align*}
   I_1 (z,F)
   &=  \int_{a-i\infty}^{a+i\infty} \hspace{-0.5cm} (2\pi Q^2)^s \left( \sum_{n=1}^\infty \frac{1}{n^{2-s}} \right) \overline{ \left( \sum_{k=1}^\infty \frac{\mu_F (k)}{k^{1-\overline{s}}} \right) } \Gamma( s - \mu )\Gamma( s + \mu )\Gamma(2-s) e^{(z - \frac{3}{2}\pi i)s}ds
   \\
   &= \sum_{k,n=1}^\infty \frac{{\mf}}{kn^2} \int_{a-i\infty}^{a+i\infty} e^{s\left\{ \log (2\pi nk Q^2) -\frac{3}{2}\pi i + z \right\}}\Gamma( s - \mu )\Gamma( s + \mu )\Gamma(2-s)ds.  \tag{4.4.6}
\end{align*}The justification of the interchange of the order of integration and summation is ensured as follows : For $0<y<3\pi$, we have 
\begin{eqnarray*}
   \lefteqn{\int_{a-i\infty}^{a+i\infty} \left| e^{s\left\{ \log (2nk\pi Q^2) + z - \frac{3}{2}\pi i \right\}}  \Gamma( s - \mu )\Gamma( s + \mu )\Gamma(2-s)  \right| |ds|}
   \\
   &&\ll (nk)^a \cdot e^{ax} \int_{-\infty}^\infty e^{-\left( y - \frac{3}{2}\pi \right)t -\frac{3}{2}\pi|t|} |t|^{a+\frac{1}{2}}dt
   \\
   &&= (nk)^a \cdot e^{ax} \left\{ \int_{-\infty}^0 e^{-\left( y - \frac{3}{2}\pi \right)t + \frac{3}{2}\pi t} (-t)^{a+\frac{1}{2}}dt + \int_{0}^\infty e^{-\left( y - \frac{3}{2}\pi \right)t - \frac{3}{2}\pi t}t^{a+\frac{1}{2}}dt  \right\}
   \\
   &&\ll_{a,x,y} (nk)^a
   \\
   &&\ll_a  (nk)^a .
\end{eqnarray*}Hence, we have
\begin{equation*}
   \sum_{k,n=1}^\infty \left| \frac{{\mf}}{kn^2} \int_{a-i\infty}^{a+i\infty} e^{s\left\{ \log (2nk\pi Q^2) + z - \frac{3}{2}\pi i \right\}}\Gamma( s - \mu )\Gamma( s + \mu )\Gamma(2-s)ds \right|
   \ll_a  \sum_{k,n=1}^\infty \left| \frac{{\mf}}{k^{1-a} n^{2-a}} \right|.
\end{equation*}Since the series for $k$ is absolutely convergent by (2.2.2), the above series of the left hand side is convergent absolutely and uniformly. Therefore, the interchange of the order of integration and summation is justified for $0 < y < 3\pi$. Let $m_1 ,m_2$ be non-negative integers. The residue of the integrand in (4.4.6) at $s = \mu - m_1$ is 
\begin{align*}
   R^{(1)}_{k,n, m_1}(z,\mu)
   &= \lim_{s \to \mu - m_1} \{ s - ( \mu - m_1 ) \}\Gamma(s-\mu)\Gamma(s+\mu)\Gamma(2-s)e^{ \left\{ \log(2\pi nkQ^2)-\frac{3}{2}\pi i + z \right\}s }
   \\
   &= \frac{(-1)^{m_1}}{m_{1}!}\Gamma(2\mu-m_1)\Gamma(2-\mu+m_1)(2\pi nkQ^2)^{\mu-m_1}e^{\left( z - \frac{3}{2}\pi i \right)( \mu - m_1 )}.  \tag{4.4.7}
\end{align*}Similarly, the residue of the integrand in (4.4.6) at $s = -\mu - m_2$ is 
\begin{equation*}
   R^{(1)}_{k,n, m_2}(z,\mu)
   = \frac{(-1)^{m_2}}{m_{2}!}\Gamma(-2\mu-m_2)\Gamma(2+\mu+m_2)(2\pi nkQ^2)^{-\mu-m_2}e^{\left( z - \frac{3}{2}\pi i \right)( -\mu - m_2 )}.  \tag{4.4.8}
\end{equation*}When $\mu=0$, $\{\Gamma(s)\}^2$ has a double pole at $s=-m$, where $m$ is a non-negative integer. For every positive $\epsilon$, using the Taylor expansion for the every factor of the integrand in (4.4.6) at $s=-m+\epsilon$, we see that the residue of the integrand in (4.4.6) at $s=-m$ is
\begin{equation*}
   R^{(1)}_{k,n,m,0}(z)
   = \frac{m+1}{m!}\frac{e^{-\left( z-\frac{3}{2}\pi i \right)m}}{(2\pi nkQ^2)^m}\left\{ \log(2\pi nkQ^2) + z -\frac{3}{2}\pi i + \sum_{k_1 = 1}^m \frac{1}{k_1} -C_0 - \frac{1}{m+1} \right\}.   \tag{4.4.9}
\end{equation*}Similarly, when $\mu=\frac{1}{2}$, $\Gamma\left( s-\frac{1}{2} \right)\Gamma\left( s+\frac{1}{2} \right) = \left( s-\frac{1}{2} \right)\left\{ \Gamma\left( s-\frac{1}{2} \right) \right\}^2$ has a double pole at $s=\frac{1}{2}-m$. In the same way, the residue of the integrand in (4.4.6) at $s=\frac{1}{2}-m$ is
\begin{align*}
   R^{(1)}_{k,n,m,\frac{1}{2}}(z)
   &= \frac{\Gamma\left( \frac{3}{2}+m \right)}{(m!)^2}\frac{e^{\left( \frac{1}{2}-m \right)\left( z -\frac{3}{2}\pi i \right)}}{(2\pi nkQ^2)^{m-\frac{1}{2}}}\left\{ m\left( \psi\left( \frac{3}{2}+m \right) -2\psi(m+1) \right)  \right.
   \\
   &- \left. m\left( \log (2\pi nkQ^2) + z -\frac{3}{2}\pi i \right) + 1 \right\},   \tag{4.4.10}
\end{align*}where $\psi(s)$ is the logarithmic derivative of $\Gamma(s)$, i.e. 
\begin{equation*}
   \psi(s) := \frac{\Gamma^\prime}{\Gamma}(s).  \tag{4.4.11}
\end{equation*}Since for $0<y<3\pi$ integrals along the upper and the lower side of the contour tend to $0$, for $\mu$ with $0 < \mu < 1$ except $\mu \neq \frac{1}{2}$, we have by theorem of residues
\footnotesize
\begin{align*}
   &\int_{a-i\infty}^{a+i\infty}e^{ s\left\{ \log (2nk\pi Q^2) + z - \frac{3}{2}\pi i \right\} } \Gamma(s - \mu )\Gamma( s + \mu )\Gamma(2-s)ds
   \\
   &= -\int_C e^{ s\left\{ \log (2nk\pi Q^2) + z - \frac{3}{2}\pi i \right\} } \Gamma(s - \mu )\Gamma( s + \mu )\Gamma(2-s)ds
   \\
   &-2\pi i\left\{ \sum_{m_1 = 0}^M R^{(1)}_{k,n, m_1}(z,\mu) + \sum_{m_2 = 0}^{M^\prime} R^{(1)}_{k,n, m_2}(z,\mu)  \right\},   \tag{4.4.12}
\end{align*}
\normalsize
where $C$ is the contour which the poles of $\Gamma(2-s)$ and those of $\Gamma( s + \mu) \Gamma ( s - \mu )$ are on opposite sides of it. Of course, when $\mu=0,\frac{1}{2}$, the terms of residue in (4.4.12) are replaced by
\begin{equation*}
   \sum_{m = 0}^{M^{\prime\prime}} R^{(1)}_{k,n,m,0}(z), \quad \sum_{m = 0}^{M^{\prime\prime}} R^{(1)}_{k,n,m,\frac{1}{2}}(z).   \tag{4.4.13} 
\end{equation*}We use the same convention hereafter. Putting $w=2-s$ in the integral round the contour $C$ on the right hand side in (4.4.12) and using Barnes type integral for the Whittaker function (4.1.6), for $\mu$ with $0 \leq \mu < 1$ and $|y-\frac{3}{2}\pi|<\frac{3}{2}\pi$, the integral round the contour $C$ in (4.4.12) is
\begin{equation*}
   2\pi i(2\pi nkQ^2)^\frac{1}{2}\exp\left( \frac{e^{\frac{3}{2}\pi i - z}}{4\pi nkQ^2} + \frac{z}{2} - \frac{3}{4}\pi i \right)\Gamma(2-\mu)\Gamma(2+\mu)
      W_{-\frac{3}{2}, \mu}\left( \frac{e^{\frac{3}{2}\pi i - z}}{2\pi nkQ^2} \right).  \tag{4.4.14}
\end{equation*}Therefore, we have
\footnotesize
\begin{align*}
   I_1 (z,F)
  &= \sum_{k,n=1}^\infty \frac{{\mf}}{kn^2} \left\{ 2\pi i \cdot (2\pi nkQ^2)^\frac{1}{2}\exp\left( \frac{e^{\frac{3}{2}\pi i - z}}{4\pi nkQ^2} + \frac{z}{2} - \frac{3}{4}\pi i \right) \right.
  \\
  &\times \left. \Gamma(2-\mu)\Gamma(2+\mu)W_{-\frac{3}{2}, \mu}\left( \frac{e^{\frac{3}{2}\pi i - z}}{2\pi nkQ^2} \right) -2\pi i \left(\sum_{m_1 = 0}^M R^{(1)}_{k,n, m_1}(z,\mu) + \sum_{m_2 = 0}^{M^\prime} R^{(1)}_{k,n, m_2}(z,\mu) \right) \right\}.   \tag{4.4.15}
\end{align*}
\normalsize
The following lemma ensures the convergence for the series on the right hand side in (4.4.15).
\begin{lemma}
For $F \in \mathcal{S}^{\text{poly}}$ with $(r,\lambda_j) = (1,1)$ and $0 \leq \mu < 1$, the series on the right hand side in (4.4.15) is absolutely and uniformly convergent on every compact subset on the whole complex plane.
\end{lemma}
We will prove Lemma 4.4.1 in the next section. By Lemma 4.4.1, for $F \in \mathcal{S}^{\text{poly}}$ whose $(r,\lambda_j) = (1,1)$ in (1.5.15) and $0 \leq \mu < 1,$ we have the following analytic continuation of $f_1 (z,F)$ for $y > -\pi$ :
\footnotesize
\begin{align*}
   f_1 (z,F) 
   &= \frac{\overline{\omega} e^{-\mu\pi i}}{(2\pi)^3 Qi}\int_a^{a+i\infty} (2\pi Q^2 )^s \frac{\zeta(2-s)}{{\f}}\Gamma(s - \mu)\Gamma(s + \mu )\Gamma(2-s)e^{\left( z + \frac{3}{2}\pi i \right)s}ds
   \\
   &-\frac{\overline{\omega} e^{\mu\pi i}}{(2\pi)^3 Qi}\int_a^{a+i\infty} (2\pi Q^2 )^s \frac{\zeta(2-s)}{{\f}}\Gamma( s - \mu )\Gamma(s + \mu)\Gamma(2-s)e^{\left(z - \frac{\pi}{2}i \right)s}ds   
   \\
   &+\frac{\overline{\omega} e^{-\mu\pi i} }{(2\pi)^3 Qi}\int_a^{a+i\infty} (2\pi Q^2 )^s \frac{\zeta(2-s)}{{\f}}\Gamma( s - \mu )\Gamma(s + \mu)\Gamma(2-s)e^{\left(z + \frac{\pi}{2}i \right)s}ds   
   \\
   &+ \frac{\overline{\omega} e^{\mu\pi i}i }{(2\pi)^3 Q} \sum_{k,n=1}^\infty \frac{{\mf}}{kn^2}\times 2\pi i\left\{(2\pi nkQ^2)^\frac{1}{2}\exp\left( \frac{e^{\frac{3}{2}\pi i - z}}{4\pi nkQ^2} + \frac{z}{2} - \frac{3}{4}\pi i \right) \right.
   \\
   &\times \left. \Gamma(2-\mu)\Gamma(2+\mu)W_{-\frac{3}{2}, \mu}\left( \frac{e^{\frac{3}{2}\pi i - z}}{2\pi nkQ^2} \right) -\sum_{m_1 = 0}^M R^{(1)}_{k,n, m_1}(z,\mu) - \sum_{m_2 = 0}^{M^\prime} R^{(1)}_{k,n, m_2}(z,\mu) \right\}
   \\
   &+ \frac{\overline{\omega} e^{\mu\pi i} }{(2\pi)^3 Qi}\int_{a-i\infty}^a (2\pi Q^2 )^s \frac{\zeta(2-s)}{{\f}}\Gamma(s - \mu )\Gamma( s + \mu )\Gamma(2-s) e^{(z - \frac{3}{2}\pi i)s}ds.   \tag{4.4.16}
\end{align*}
\normalsize
The first is analytic for $y>-3\pi$, the second for $y>-\pi$, the third for $y>-2\pi$, the fourth is analytic on the whole complex plane by Lemma 4.4.1, and the next is analytic for $y<3\pi$. Therefore, (4.4.16) completes the proof of the continuation of $f(z,F)$ to the region $y>-\pi$.   $\square$ 

\subsection{Proof of Lemma 4.4.1}
We prove Lemma 4.4.1. First, we consider the case $\mu > \frac{1}{2}$. By the asymptotic expansion (4.1.12),
\begin{equation*}
   W_{-\frac{3}{2}, \mu}\left( \frac{e^{\frac{3}{2}\pi i - z}}{2\pi nkQ^2} \right)
   = \frac{\Gamma(2\mu)}{\Gamma(2+\mu)} \left( \frac{e^{\frac{3}{2}\pi i - z}}{2\pi nkQ^2} \right)^{\frac{1}{2}-\mu} 
   + O\left( \frac{|e^{ \left( \frac{3}{2}\pi i - z \right)\left( \frac{3}{2} -\mu \right) }|}{(2\pi nkQ^2)^{\frac{3}{2} - \mu}} \right).
\end{equation*}Hence, the inside of the curly brackets on the right hand side of (4.4.15) is
\begin{eqnarray*}
   \lefteqn{2\pi i(2\pi nkQ^2)^\frac{1}{2} \exp\left( \frac{e^{\frac{3}{2}\pi i - z}}{4\pi nkQ^2} +\frac{z}{2} - \frac{3}{4}\pi i \right)\Gamma(2-\mu)\Gamma(2+\mu)W_{-\frac{3}{2}, \mu}\left( \frac{e^{\frac{3}{2}\pi i - z}}{2\pi nkQ^2} \right)}
    \\
    &&-2\pi i\left\{ \sum_{m_1 =0}^M R^{(1)}_{k,n, m_1}(z,\mu) + \sum_{m_2 = 0}^{M^\prime} R^{(1)}_{k,n, m_2}(z,\mu)  \right\}
    \\
    &&= 2\pi i(2\pi nkQ^2)^\mu \exp\left( \frac{e^{\frac{3}{2}\pi i - z}}{4\pi nkQ^2} -\frac{3}{2}\mu\pi i + \mu z \right)\Gamma( 2 - \mu)\Gamma( 2\mu)
    \\
    &&+O_{Q,\mu,x} \left( \frac{1}{(2\pi nkQ^2)^{1 - \mu}} \exp\left(\frac{e^{-x}}{4\pi nkQ^2} \right) \right) 
    \\
    &&- 2\pi i\left\{ \sum_{m_1 = 0}^M R^{(1)}_{k,n, m_1}(z,\mu) + \sum_{m_2 = 0}^{M^\prime} R^{(1)}_{k,n, m_2}(z,\mu)  \right\}. 
\end{eqnarray*}Now, by the Taylor expansion
\begin{align*}
   \exp\left( \frac{e^{\frac{3}{2}\pi i - z}}{4\pi nkQ^2} \right) 
   &= 1 + O\left( \left| \frac{e^{\frac{3}{2}\pi i - z}}{4\pi nkQ^2} \right| \right)
   \\
   &= 1 + O\left( \frac{e^{-x}}{4\pi nkQ^2} \right)   \tag{4.5.1}
\end{align*}as $n,k$ tend to infinity and (4.4.7), (4.4.8), we have
\begin{eqnarray*}
    \lefteqn{2\pi i(2\pi nkQ^2)^\frac{1}{2} \exp\left( \frac{e^{\frac{3}{2}\pi i - z}}{4\pi nkQ^2} +\frac{z}{2} - \frac{3}{4}\pi i \right)\Gamma(2-\mu)\Gamma(2+\mu)W_{-\frac{3}{2}, \mu}\left( \frac{e^{\frac{3}{2}\pi i - z}}{2\pi nkQ^2} \right)}
    \\
    &&-2\pi i\left\{ \sum_{m_1 = 0}^M R^{(1)}_{k,n, m_1}(z,\mu) + \sum_{m_2 = 0}^{M^\prime} R^{(1)}_{k,n, m_2}(z,\mu)  \right\}
    \\
    &&= O_{Q,\mu, x}\left( \frac{1}{(nk)^{1 - \mu}} \right) + O_{Q,\mu,x}\left( \frac{1}{(nk)^{1 - \mu}} + \frac{1}{(nk)^{2 - \mu}} \right)
    \\
    &&-2\pi i\left\{ \sum_{m_1 = 1}^M R^{(1)}_{k,n, m_1}(z,\mu) + \sum_{m_2 = 0}^{M^\prime} R^{(1)}_{k,n, m_2}(z,\mu)  \right\}
    \\
    &&\ll O_{Q,\mu, x}\left( \frac{1}{(nk)^{1 - \mu}} \right) + \left\{ \sum_{m_1 = 1}^M \left| R^{(1)}_{k,n, m_1}(z,\mu) \right| + \sum_{m_2 = 0}^{M^\prime} \left| R^{(1)}_{k,n, m_2}(z,\mu) \right|  \right\}
    \\
    &&\ll_{Q,\mu,x} O_{Q,\mu, x}\left( \frac{1}{(nk)^{1 - \mu}} \right) + \sum_{m_1 = 1}^M \frac{1}{(nk)^{m_1 - \mu}} + \sum_{m_2 = 0}^{M^\prime} \frac{1}{(nk)^{m_2 + \mu}}
    \\
    &&= O_{Q,\mu,x}\left( \frac{1}{(nk)^{1 - \mu}} \right)
\end{eqnarray*}Hence, $I_1 (z,F)$ is evaluated as follows :
\begin{align*}
   I_1 (z,F)
   &\ll_{Q,\mu,x} \sum_{k,n=1}^\infty \frac{|{\mf}|}{kn^2} \cdot O_{Q,\mu, x}\left( \frac{1}{(nk)^{1 - \mu}} \right)   
   \\
   &= O_{Q,\mu,x}\left(\sum_{k,n=1}^\infty \frac{|{\mf}|}{k^{2 - \mu}n^{3 - \mu}} \right).
\end{align*}Therefore, the series on the right hand side in (4.4.15) is convergent for $\frac{1}{2} < \mu < 1$. 

Secondly, in the case $\mu = \frac{1}{2}$, by (4.1.13),
\begin{align*}
  W_{-\frac{3}{2},\frac{1}{2}}\left( \frac{e^{\frac{3}{2}\pi i - z}}{2\pi nkQ^2} \right)
   &= \frac{1}{\Gamma\left( \frac{5}{2} \right)} + O_y \left( \frac{e^{-x}}{2\pi nkQ^2}\log \left( \frac{e^{-x}}{2\pi nkQ^2} \right) \right)
   \\
   &= \frac{1}{\Gamma\left( \frac{5}{2} \right)} + O_y \left( \left( \frac{e^{-x}}{2\pi nkQ^2} \right)^{1-\delta} \right),
\end{align*}where $\delta$ is any positive real number. By the same calculation as in the first case and using (4.4.10), the inside of the curly brackets on the right hand side of (4.4.15) is evaluated as follows :
\begin{eqnarray*}
   \lefteqn{2\pi i(2\pi nkQ^2)^\frac{1}{2} \exp\left( \frac{e^{\frac{3}{2}\pi i - z}}{4\pi nkQ^2} +\frac{z}{2} - \frac{3}{4}\pi i \right)\Gamma\left( \frac{3}{2} \right)\Gamma\left( \frac{5}{2} \right)W_{-\frac{3}{2}, \frac{1}{2}}\left( \frac{e^{\frac{3}{2}\pi i - z}}{2\pi nkQ^2} \right)}
   \\
   &&- 2\pi i\sum_{m= 0}^{M^{\prime\prime}} R^{(1)}_{k,n, m,\frac{1}{2}}(z)
   \\
    &&\ll O_{Q, x,y}\left( \frac{1}{(nk)^{\frac{1}{2}-\delta}} \right) + \sum_{m=1}^{M^{\prime\prime}} \frac{\log nk}{(nk)^{m - \frac{1}{2}}} 
    \\
    &&= O_{Q,x,y,M^{\prime\prime}}\left(\frac{1}{(nk)^{\frac{1}{2}-\delta}}\right).
\end{eqnarray*}Hence, $I_1 (z,F)$ is evaluated as follows :
\begin{align*}
I_1 (z,F)
&\ll \sum_{k,n=1}^\infty \frac{|{\mf}|}{kn^2} \cdot O_{Q,x,y,M^{\prime\prime}}\left( \frac{1}{(nk)^{\frac{1}{2}-\delta}}\right)
\\
&= O_{Q,x,y,M^{\prime\prime}} \left( \sum_{k,n=1}^\infty \frac{|{\mf}|}{k^{\frac{3}{2}-\delta}n^{\frac{5}{2}-\delta}} \right).
\end{align*}
Therefore, the series on the right hand side in (4.4.15) is convergent for $\mu = \frac{1}{2}$. 

Thirdly, in the case $0 < \mu < \frac{1}{2}$, by (4.1.14),
\begin{equation*}
   W_{-\frac{3}{2}, \mu}\left( \frac{e^{\frac{3}{2}\pi i - z}}{2\pi nkQ^2} \right)
   = \frac{\Gamma(2\mu)}{\Gamma(2+\mu)} \left( \frac{e^{\frac{3}{2}\pi i - z}}{2\pi nkQ^2} \right)^{\frac{1}{2}-\mu} 
   + O\left( \frac{e^{-x\left( \mu + \frac{1}{2} \right)}}{(2\pi nkQ^2)^{ \mu + \frac{1}{2}}} \right)
\end{equation*}and
\begin{eqnarray*}
   \lefteqn{2\pi i(2\pi nkQ^2)^\frac{1}{2} \exp\left( \frac{e^{\frac{3}{2}\pi i - z}}{4\pi nkQ^2} +\frac{z}{2} - \frac{3}{4}\pi i \right)\Gamma(2-\mu)\Gamma(2+\mu)W_{-\frac{3}{2}, \mu}\left( \frac{e^{\frac{3}{2}\pi i - z}}{2\pi nkQ^2} \right)}
    \\
    &&-2\pi i\left\{ \sum_{m_1 = 0}^M R^{(1)}_{k,n, m_1}(z,\mu) + \sum_{m_2 = 0}^{M^\prime} R^{(1)}_{k,n, m_2}(z,\mu)  \right\}
    \\
    &&\ll_{Q,\mu,x} O_{Q,\mu, x}\left( \frac{1}{(nk)^{1 - \mu}} \right) + \sum_{m_1 = 1}^M \frac{1}{(nk)^{m_1 - \mu}} + \sum_{m_2 = 0}^{M^\prime} \frac{1}{(nk)^{m_2 + \mu}}
    \\
    &&= O_{Q,\mu, x}\left( \frac{1}{(nk)^{1 - \mu}} \right).
\end{eqnarray*}Therefore, the series on the right hand side in (4.4.15) is convergent for $0 < \mu < \frac{1}{2}$.

Finally, in the case $\mu=0$, by (4.1.16),
\begin{align*}
   W_{-\frac{3}{2},0}\left( \frac{e^{\frac{3}{2}\pi i - z}}{2\pi nkQ^2} \right)
   &= -\frac{e^{\frac{3}{4}\pi i -\frac{z}{2}}}{(2\pi nkQ^2)^\frac{1}{2}}\left\{ \log\left( \frac{e^{-x}}{2\pi nkQ^2} \right) + \left( \frac{3}{2}\pi - y \right)i + \frac{\Gamma^\prime}{\Gamma}(2) + 2C_0 \right\}
   \\
   &+ O_{Q,x} \left( \frac{1}{(nk)^\frac{3}{2}}\log\left( \frac{e^{-x}}{2\pi nkQ^2} \right) \right).
\end{align*}Using the recurrence formula
\begin{equation*}
   \psi(s+1) = \frac{1}{s} + \psi(s)   \tag{4.5.2}
\end{equation*}(see~\cite{L}) and (4.1.16), the inside of the curly brackets on the right hand side of (4.1.15) is evaluated as follows :
\begin{eqnarray*}
   \lefteqn{2\pi i (2\pi nkQ^2)^\frac{1}{2} \exp\left( \frac{e^{\frac{3}{2}\pi i - z} }{4\pi nkQ^2} + \frac{z}{2} - \frac{3}{4}\pi i \right) \Gamma(2)^2 W_{ -\frac{3}{2},0 } \left( \frac{e^{\frac{3}{2}\pi i - z} }{2\pi nkQ^2} \right) - 2\pi i\sum_{m = 0}^{ M^{\prime\prime} } R^{(1)}_{k,n, m,0}(z)}
   \\
   &&= -2\pi i \left( \frac{\Gamma^\prime}{\Gamma}(2) + C_0 -1 \right) + O_{Q,x,y} \left( \frac{1}{nk}\log \left( \frac{e^{-x}}{2\pi nkQ^2} \right) \right) -2\pi i\sum_{m = 1}^{ M^{\prime\prime} } R^{(1)}_{k,n,m,0}(z)
    \\
    &&= O_{Q,x,y} \left( \frac{1}{nk}\log \left( \frac{e^{-x}}{2\pi nkQ^2} \right) \right) - 2\pi i\sum_{m = 1}^{ M^{\prime\prime} } R^{(1)}_{k,n,m,0}(z)
    \\
    &&= O_{Q,x,y} \left( \frac{1}{nk} \cdot \frac{1}{(nk)^{1-\delta}} \right) - 2\pi i\sum_{m = 1}^{ M^{\prime\prime} } R^{(1)}_{k,n,m,0}(z)
    \\
    &&= O_{Q,x,y}\left( \frac{1}{(nk)^{2-\delta}} \right) -2\pi i\sum_{m = 1}^{M^{\prime\prime}} R^{(1)}_{k,n,m,0}(z)
    \\
    &&\ll_{Q,x,y}  O_{Q,x,y} \left( \frac{1}{(nk)^{2-\delta}} \right) + \sum_{m = 1}^{M^{\prime\prime}} \frac{\log nk}{(nk)^m}
    \\
    &&\ll O_{M^{\prime\prime}}\left( \frac{1}{(nk)^{1-\delta}} \right),
\end{eqnarray*}where $\delta$ is any positive real number. Hence, $I_1 (z,F)$ is evaluated as follows :
\begin{align*}
   I_1 (z,F)
   &\ll_{M^{\prime\prime}}  \sum_{k,n=1}^\infty \frac{|{\mf}|}{kn^2} \cdot \frac{1}{(nk)^{1-\delta}}
   \\
   &= \sum_{k,n=1}^\infty \frac{|{\mf}|}{k^{2-\delta}n^{3-\delta}}.
\end{align*}Therefore, the series on the right hand side in (4.4.15) is convergent for $\mu = 0$.  In summary, 
\begin{align*}
   I_1 (z,F) \leq \sum_{k,n=1}^\infty \frac{|{\mf}|}{kn^2} \cdot \max\left\{ \frac{1}{(nk)^{ 1 - \mu }},\frac{1}{(nk)^{\frac{1}{2}-\delta}} \right\}.  \tag{4.5.3}
\end{align*}By (4.1.2),(4.1.4) and (4.1.5), the Whittaker function $W_{-\frac{3}{2},\mu} \left( \frac{e^{\frac{3}{2}\pi i-z} }{2\pi nkQ^2} \right)$ is analytic for all $z \in \mathbb{C}$. Therefore, We have the desired result. $\square$ 

\subsection{Two proofs of Corollary 2.2.3}
We prove Corollary 2.2.3. We use the following lemma similar to Lemma 4.2.1. We can prove this lemma by modifying the proof of Lemma 4.2.1.
\begin{lemma}
Let $F \in \mathcal{S}$ and let $T$ be sufficiently large.Moreover, let $H=D\log\log T$ be fixed, where $D$ is a large positive constant.In any subinterval of length $1$ in $[-T-H, -T+H]$ there are lines $t=t_0$ such that
\begin{equation*}
   |F( \sigma+it_0 )|^{-1} = O(\exp(C(\log T)^2))   \tag{4.6.1}  
\end{equation*}uniformly in $\sigma \geq -2$.
\end{lemma}
We consider the integral
\begin{equation*}
   \int_{\mathscr{L}^\prime} \frac{\zeta(s-1)}{F(s)}e^{zs}ds,  \tag{4.6.2}
\end{equation*}where $\mathscr{L}^\prime$ is the contour symmetrical upon the real axis to $\mathscr{L}$ in (4.2.3). By Lemma 4.6.1, the integral along the lower side of the contour tends to $0$ as $n$ tends to infinity for $z \in \mathbb{H}^{-}$. Then, we have by residue theorem and the definition (2.2.9), in a similar manner as (4.2.4),
\begin{equation*}
   2\pi i f^{-}(z,F) = f_{1}^{-}(z,F) +  f_{2}^{-}(z,F) +  f_{3}^{-}(z,F),   \tag{4.6.3}
\end{equation*}where
\begin{equation*}
   f_{1}^{-}(z,F) = \int_{a}^{a-i\infty} \frac{\zeta(s-1)}{F(s)}e^{sz}ds   \tag{4.6.4} 
\end{equation*}is analytic on $\mathbb{H}^-$, 
\begin{equation*}
  f_{2}^{-}(z,F) = \int_{\overline{L}} \frac{\zeta(s-1)}{F(s)}e^{sz}ds   \tag{4.6.5}
\end{equation*}is analytic on the whole complex plane. We consider the same setting as in $L$ for the curve $\overline{L}$. In the same way as obtaining (4.3.1),
\begin{align*}
   f_{3}^{-}(z,F)
   &= \int_{b - i\infty}^b \left( \sum_{n=1}^\infty \frac{g(n)}{n^s} \right)e^{zs}ds
   \\
   &= \sum_{n=1}^\infty g(n) \int_{b - i\infty}^b e^{s(z-\log n)}ds
   \\
   &= e^{bz}\sum_{n=1}^\infty \frac{g(n)}{n^b (z-\log n)}   \tag{4.6.6}
\end{align*}is meromorphic on the whole complex plane. Now we already know that $f_{1}^{-} (z,F)$ is analytic for $y<0$, and we have to continue to $y<\pi$ just as in the case of $f_1 (z,F)$ (see Section 4.4). By the functional equation for $\zeta(s)$ and $F(s)$, we have
\begin{align*}
  f_{1}^{-} (z,F)
   &= \int_{a}^{a-i\infty} \frac{\zeta(s-1)}{F(s)}e^{zs}ds
   \\
   &= -\int_{a-i\infty}^a \frac{\zeta(s-1)}{F(s)}e^{zs}ds
   \\
   &= f_{11}^- (z,F) + f_{12}^- (z,F) + f_{13}^- (z,F) + f_{14}^- (z,F),  \tag{4.6.7}
\end{align*}where
\begin{equation*}
   f_{11}^{-}(z,F) 
   = \frac{\overline{\omega} e^{-\mu\pi i} }{(2\pi)^3 Qi}\int_{a-i\infty}^a (2\pi Q^2 )^s \frac{\zeta(2-s)}{{\f}}\Gamma(s+\mu)\Gamma(s-\mu)\Gamma(2-s)e^{\left( z + \frac{3}{2}\pi i \right)s}ds   \tag{4.6.8}
\end{equation*}is analytic for $y < 0$,
\begin{equation*}
   f_{12}^{-}(z,F)
   = -\frac{\overline{\omega} e^{\mu\pi i} }{(2\pi)^3 Qi}\int_{a-i\infty}^a (2\pi Q^2 )^s \frac{\zeta(2-s)}{{\f}}\Gamma(s+\mu)\Gamma(s-\mu)\Gamma(2-s)e^{\left( z-\frac{\pi}{2}i \right)s}ds    \tag{4.6.9}
\end{equation*}for $y < 2\pi$, 
\begin{equation*}
   f_{13}^{-}(z,F)
  =  \frac{\overline{\omega} e^{-\mu\pi i} }{(2\pi)^3 Qi} \int_{a-i\infty}^a (2\pi Q^2 )^s \frac{\zeta(2-s)}{{\f}}\Gamma(s+\mu)\Gamma(s-\mu)\Gamma(2-s)e^{\left( z+\frac{\pi}{2}i \right)s}ds    \tag{4.6.10}
\end{equation*}for $y < \pi$,
\begin{equation*}
  f_{14}^{-}(z,F)
  = -\frac{\overline{\omega} e^{\mu\pi i}}{(2\pi)^3 Qi} \int_{a-i\infty}^a (2\pi Q^2 )^s \frac{\zeta(2-s)}{{\f}}\Gamma(s+\mu)\Gamma(s-\mu)\Gamma(2-s)e^{\left( z-\frac{3}{2}\pi i \right)s}ds   \tag{4.6.11}
\end{equation*} for $y < 3\pi$. Splitting the integral on the right hand side in (4.6.8) just as in the case of $f_{14}(z,F)$, we have
\[ f_{11}^{-}(z,F) = I_{1}^{-} (z,F) + I_{2}^{-} (z,F), \]where
\begin{equation*}
    I_{1}^{-} (z,F) = \frac{\overline{\omega}e^{-\mu\pi i}}{(2\pi)^3 Qi} \int_{a-i\infty}^{a+i\infty}(2\pi Q^2)^s \frac{\zeta(2-s)}{{\f}}\Gamma(s+\mu)\Gamma(s-\mu)\Gamma(2-s)e^{\left( z + \frac{3}{2}\pi i \right)s}ds   \tag{4.6.12}
\end{equation*}and
\begin{equation*}
    I_{2}^{-} (z,F) = -\frac{\overline{\omega}e^{-\mu\pi i}}{(2\pi)^3 Qi} \int_{a}^{a+i\infty}(2\pi Q^2)^s \frac{\zeta(2-s)}{{\f}}\Gamma(s+\mu)\Gamma(s-\mu)\Gamma(2-s)e^{\left( z + \frac{3}{2}\pi i \right)s}ds.   \tag{4.6.13}
 \end{equation*}We see that the integral $I_{2}^{-} (z,F)$ is convergent for $y > -3\pi$ by the same way as in (4.4.5). Since $f_{11}^{-}(z,F)$ is analytic for $y < 0$, the integral $I_{1}^{-} (z,F)$ is convergent for $-3\pi < y < 0$ and we can calculate $I_{1}^{-}(z,F)$ for $-3\pi < y < 0$ in a similar way as (4.4.4) (Section 4.4). Let $m_1, m_2$ and $m$ be non-negative integers. By taking the path of integration $C$ in (4.4.12), we have for $0 \leq \mu < 1$ and $|y+\frac{3}{2}\pi| < \frac{3}{2}\pi$
\footnotesize
\begin{align*}
    I_{1}^{-}(z,F) 
    &= \frac{\overline{\omega}e^{-\mu\pi i}}{(2\pi)^3 Qi} \sum_{k,n=1}^\infty \frac{{\mf}}{kn^2} \times 2\pi i \left\{(2\pi nkQ^2)^\frac{1}{2}\exp\left( \frac{3}{4}\pi i + \frac{z}{2} + \frac{e^{-\frac{3}{2}\pi i-z} }{4\pi nkQ^2} \right) \right.
    \\
    &\times \left. \Gamma(2+\mu)\Gamma(2-\mu)W_{-\frac{3}{2}, \> \mu}\left( \frac{e^{-\frac{3}{2}\pi i - z}}{2\pi nkQ^2} \right) - \sum_{m_1 = 0}^M R^{(2)}_{k,n, m_1}(z,\mu) \hspace{-0.01cm} - \hspace{-0.01cm} \sum_{m_2 = 0}^{M^\prime} R^{(2)}_{k,n, m_2}(z,\mu) \hspace{-0.01cm} \right\},  \tag{4.6.14}
\end{align*}
\normalsize 
where
\begin{gather*}
   R^{(2)}_{k,n, m_1}(z,\mu)
    = \frac{(-1)^{m_1}}{m_{1}!}\Gamma(2\mu-m_1)\Gamma(2-\mu+m_1)(2\pi nkQ^2)^{\mu-m_1}e^{\left( z + \frac{3}{2}\pi i \right)( \mu - m_1 )},  \tag{4.6.15}
    \\
    R^{(2)}_{k,n, m_2}(z,\mu)
    = \frac{(-1)^{m_2}}{m_{2}!}\Gamma(-2\mu-m_2)\Gamma(2+\mu+m_2)(2\pi nkQ^2)^{-\mu-m_2}e^{\left( z + \frac{3}{2}\pi i \right)( -\mu - m_2 )},  \tag{4.6.16}
    \\
    R^{(2)}_{k,n,m,0}(z)
    =  \frac{m+1}{m!}\frac{e^{-\left( z+\frac{3}{2}\pi i \right)m}}{(2\pi nkQ^2)^m}\left\{ \log(2\pi nkQ^2) + z  + \frac{3}{2}\pi i  + \sum_{k_1 = 1}^m \frac{1}{k_1} - C_0 - \frac{1}{m+1} \right\}  \tag{4.6.17}
\end{gather*}and
\begin{align*}
   R^{(2)}_{k,n,m,\frac{1}{2}}(z)
   &= \frac{\Gamma\left( \frac{3}{2}+m \right)}{(m!)^2}\frac{e^{\left( \frac{1}{2}-m \right)\left( z + \frac{3}{2}\pi i \right)}}{(2\pi nkQ^2)^{m-\frac{1}{2}}}\left\{ m\left( \psi\left( \frac{3}{2}+m \right) -2\psi(m+1) \right)  \right.
   \\
   &- \left. m\left( \log (2\pi nkQ^2) + z + \frac{3}{2}\pi i \right) + 1 \right\}      \tag{4.6.18}
\end{align*}
 are residues of the integrand in (4.6.12) at  $s = \mu - m_1, -\mu - m_2, -m$ and $\frac{1}{2}-m$ respectively. The convergence of the series on the right hand side in (4.6.14) follows in a similar manner as the consideration in (4.4.15). Finally, by (4.6.7)-(4.6.18) we obtain the following continuation of $f_{1}^{-}(z,F)$ to $y<\pi$ : For $F \in \mathcal{S}^{\text{poly}}$ whose $(r,\lambda_j) = (1,1)$ in (1.5.15) and $0 \leq \mu < 1$, 
\footnotesize
\begin{align*}
    f_{1}^{-} (z,F)
    &= \frac{\overline{\omega}e^{-\mu\pi i}}{(2\pi)^3 Qi} \sum_{k,n=1}^\infty \frac{{\mf}}{kn^2}  \times 2\pi i \left\{ (2\pi nkQ^2)^\frac{1}{2}\exp\left( \frac{3}{4}\pi i + \frac{z}{2} + \frac{e^{-\frac{3}{2}\pi i-z} }{4\pi nkQ^2} \right)  \right.
    \\
    &\times \left. \Gamma(2+\mu)\Gamma(2-\mu)W_{-\frac{3}{2}, \> \mu}\left( \frac{e^{-\frac{3}{2}\pi i - z}}{2\pi nkQ^2} \right) - \sum_{m_1 = 0}^M R^{(2)}_{k,n, m_1}(z,\mu) - \sum_{m_2 = 0}^{M^\prime} R^{(2)}_{k,n, m_2}(z,\mu) \right\}
    \\
    &- \frac{\overline{\omega}e^{-\overline{\mu}\pi i}}{(2\pi)^3 Qi}\int_{a}^{a+i\infty} (2\pi Q^2)^s \frac{\zeta(2-s)}{{\f}}\Gamma(s+\mu)\Gamma(s-\mu)\Gamma(2-s)e^{\left( z + \frac{3}{2}\pi i \right)s}ds
    \\
    &-\frac{\overline{\omega} e^{\mu\pi i} }{(2\pi)^3 Qi}\int_{a-i\infty}^a (2\pi Q^2 )^s \frac{\zeta(2-s)}{{\f}}\Gamma(s+\mu)\Gamma(s-\mu)\Gamma(2-s)e^{\left( z-\frac{\pi}{2}i \right)s}ds
    \\
    &+ \frac{\overline{\omega} e^{-\mu\pi i} }{(2\pi)^3 Qi} \int_{a-i\infty}^a (2\pi Q^2 )^s \frac{\zeta(2-s)}{{\f}}\Gamma(s+\mu)\Gamma(s-\mu)\Gamma(2-s)e^{\left( z+\frac{\pi}{2}i \right)s}ds
    \\
    &-\frac{\overline{\omega} e^{\mu\pi i}}{(2\pi)^3 Qi} \int_{a-i\infty}^a (2\pi Q^2 )^s \frac{\zeta(2-s)}{{\f}}\Gamma(s+\mu)\Gamma(s-\mu)\Gamma(2-s)e^{\left( z-\frac{3}{2}\pi i \right)s}ds.  \tag{4.6.19}
\end{align*}
\normalsize
Since a lemma similar to Lemma 4.4.1 holds for the series on the right hand side in (4.6.19), the series we now consider is also absolutely and uniformly convergent on every compact subset on the whole complex plane. Therefore, we complete the continuation of $f^{-}(z,F)$ analytic for $y<0$ to the region $y<\pi$.

Also, Corollary 2.2.3 can be proved form Theorem 2.2.2 and the definition (2.2.9) directly as follows : For $z \in \mathbb{H}^-$ and $\rho$\hspace{0.01cm} with $\text{Im} \hspace{0.1cm} \rho <0$, 
\begin{align*}
   f^- (z,F)
   &= \sum_{\rho} \frac{e^{\rho z}\zeta(\rho-1)}{F^\prime (\rho)} 
   \\
   &=  \overline{\sum_{\rho} \frac{ \overline{\zeta(\rho-1)} }{\overline{F^\prime (\rho)} }e^{\overline{\rho z} } }
   \\
   &=  \overline{\sum_{\rho^\prime } \frac{\zeta(\rho^\prime -1)}{\overline{F^\prime \left(\overline{\rho^\prime } \right) } }e^{\rho^\prime \overline{z} } },
\end{align*}where $\rho^\prime = \overline{\rho}$. We recall the definition $\overline{F} (s) = \overline{F(\overline{s})}$ for $F \in {\se}$. Hence, the sum on the right hand side in the third line yields
\begin{equation*}
   \overline{\sum_{\rho^\prime} \frac{\zeta(\rho^\prime -1)}{\overline{F^\prime}(\rho^\prime)  }e^{\rho^\prime \overline{z} } }
   = f(\overline{z},\overline{F}).
\end{equation*}Here, we use the fact that if $F \in {\s}$, then so is $\overline{F} \in {\s}$. Of course, ${\s}$ may be replaced by ${\se}$. By Theorem 2.2.2, the function $f(z,F)$ has a meromorphic continuation to $y > -\pi$. Hence $f(\overline{z},\overline{F})$ has a meromorphic continuation to $y < \pi$.

\subsection{Proof of Theorem 2.2.4}
We assume the condition written in the statement of Theorem 2.2.4, that is, $F \in \mathcal{S}^{\text{poly}}$ whose $(r,\lambda_j) = (1,1)$ in (1.5.15) and $0 \leq \mu < 1$.We add (4.4.16) to (4.6.19). Since some integrals are canceled, we have for $|y| < \pi$
\footnotesize
\begin{align*}
   f_1 (z,F) + f_{1}^{-}(z,F)
   &= f_{11}(z,F) + f_{12}(z,F) + f_{13}(z,F) + \frac{\overline{\omega}e^{\mu\pi i}}{(2\pi)^3 Qi}I_1 (z,F) +  \frac{\overline{\omega}e^{\mu\pi i}}{(2\pi)^3 Qi}I_2 (z,F) 
   \\
   &+ I_{1}^- (z,F) + I_{2}^- (z,F) + f_{12}^- (z,F) + f_{13}^- (z,F) + f_{14}^- (z,F)  
   \\ 
   &= \frac{\overline{\omega}e^{\mu\pi i}i}{(2\pi)^3 Q} \sum_{k,n=1}^\infty \frac{{\mf}}{kn^2}  \times 2\pi i \left\{ (2\pi nkQ^2)^\frac{1}{2}\exp\left( -\frac{3}{4}\pi i + \frac{z}{2} + \frac{e^{\frac{3}{2}\pi i-z} }{4\pi nkQ^2} \right)  \right.
   \\
   &\times \left. \Gamma(2+\mu)\Gamma(2-\mu)W_{-\frac{3}{2}, \> \mu}\left( \frac{e^{\frac{3}{2}\pi i - z}}{2\pi nkQ^2} \right) - \sum_{m_1 = 0}^M R^{(1)}_{k,n, m_1}(z,\mu) - \sum_{m_2 = 0}^{M^\prime} R^{(1)}_{k,n, m_2}(z,\mu) \right\}
  \\
  &+\frac{\overline{\omega}e^{-\mu\pi i}}{(2\pi)^3 Qi} \sum_{k,n=1}^\infty \frac{{\mf}}{kn^2} \times 2\pi i \left\{ (2\pi nkQ^2)^\frac{1}{2} \exp\left( \frac{3}{4}\pi i + \frac{z}{2} + \frac{e^{-\frac{3}{2}\pi i-z} }{4\pi nkQ^2} \right) \right.
  \\
  &\times \left. \Gamma(2+\mu)\Gamma(2-\mu) W_{-\frac{3}{2}, \> \mu}\left( \frac{e^{-\frac{3}{2}\pi i - z}}{2\pi nkQ^2} \right) - \sum_{m_1 = 0}^M R^{(2)}_{k,n, m_1}(z,\mu) - \sum_{m_2 = 0}^{M^\prime} R^{(2)}_{k,n, m_2}(z,\mu) \right\}
   \\
   &+ A_{1}(z,F) + A_{2}(z,F),
\end{align*}
\normalsize
where
\begin{equation*}
   A_{1} (z,F) = -\frac{\overline{\omega} e^{\mu\pi i}}{(2\pi)^3 Qi} \int_{a-i\infty}^{a+i\infty} (2\pi Q^2 )^s \frac{\zeta(2-s)}{{\f}}\Gamma(s+\mu)\Gamma(s-\mu)\Gamma(2-s)e^{\left( z - \frac{\pi}{2}i \right)s}ds   \tag{4.7.1}   
\end{equation*}and
\begin{equation*}
   A_{2} (z,F) = \frac{\overline{\omega} e^{-\mu\pi i}}{(2\pi)^3 Qi} \int_{a-i\infty}^{a+i\infty} (2\pi Q^2 )^s \frac{\zeta(2-s)}{{\f}}\Gamma(s+\mu)\Gamma(s-\mu)\Gamma(2-s)e^{\left( z + \frac{\pi}{2}i \right)s}ds.  \tag{4.7.2}
\end{equation*}The integrals $A_1 (z,F)$ and $A_2 (z,F)$ are convergent for $|y| < \pi$ and we can obtain series expressions of them involving Whittaker functions in a way similar to the case of $I_1 (z,F)$ and $I_{1}^{-}(z,F)$. We have for $|y|<\pi, F \in \mathcal{S}^{\text{poly}}$ whose $(r,\lambda_j) = (1,1)$ in (1.5.15) and $0 \leq \mu < 1$
\begin{align*}
    A_1 (z,F)
    &= -\frac{\overline{\omega}e^{\mu\pi i}}{(2\pi)^3 Qi} \sum_{k,n=1}^\infty \frac{{\mf}}{kn^2} \times 2\pi i
    \left\{ (2\pi nkQ^2)^\frac{1}{2}\exp\left( -\frac{\pi}{4}i + \frac{z}{2} + \frac{e^{\frac{\pi}{2}i-z} }{4\pi nkQ^2} \right)  \right.
    \\
    &\times \left. \Gamma(2+\mu)\Gamma(2-\mu)W_{-\frac{3}{2}, \> \mu} \left( \frac{e^{\frac{\pi}{2}i - z}}{2\pi nkQ^2} \right) -\sum_{m_1 = 0}^M R^{(3)}_{k,n, m_1}(z,\mu) - \sum_{m_2 = 0}^{M^\prime} R^{(3)}_{k,n, m_2}(z,\mu) \right\},
\end{align*}
\normalsize
where
\begin{gather*}
   R^{(3)}_{k,n, m_1}(z,\mu)
    = \frac{(-1)^{m_1}}{m_{1}!}\Gamma(2\mu-m_1)\Gamma(2-\mu+m_1)(2\pi nkQ^2)^{\mu-m_1}e^{\left( z - \frac{\pi}{2}i \right)( \mu - m_1 )},  \tag{4.7.3}
    \\
    R^{(3)}_{k,n, m_2}(z,\mu)
    = \frac{(-1)^{m_2}}{m_{2}!}\Gamma(-2\mu-m_2)\Gamma(2+\mu+m_2)(2\pi nkQ^2)^{-\mu-m_2}e^{\left( z - \frac{\pi}{2}i \right)( -\mu - m_2 )},  \tag{4.7.4}
    \\
    R^{(3)}_{k,n,m,0}(z)
    =  \frac{m+1}{m!}\frac{e^{-\left( z-\frac{\pi}{2}i \right)m}}{(2\pi nkQ^2)^m}\left\{ \log(2\pi nkQ^2) + z - \frac{\pi}{2}i + \sum_{k_1 = 1}^m \frac{1}{k_1} - C_0 - \frac{1}{m+1} \right\}   \tag{4.7.5}
\end{gather*}and
\begin{align*}
   R^{(3)}_{k,n,m,\frac{1}{2}}(z)
   &= \frac{\Gamma\left( \frac{3}{2}+m \right)}{(m!)^2}\frac{e^{\left( \frac{1}{2}-m \right)\left( z -\frac{\pi}{2} i \right)}}{(2\pi nkQ^2)^{m-\frac{1}{2}}}\left\{ m\left( \psi\left( \frac{3}{2}+m \right) -2\psi(m+1) \right)  \right.
   \\
   &- \left. m\left( \log (2\pi nkQ^2) + z - \frac{\pi}{2}i \right) + 1 \right\} \tag{4.7.6}
\end{align*}are residues of the integrand in $A_1(z,F)$ at  $s = \mu - m_1, -\mu - m_2, -m$ and $\frac{1}{2}-m$ respectively. Similarly, for $|y|<\pi, F \in \mathcal{S}^{\text{poly}}$ whose $(r,\lambda_j) = (1,1)$ in (1.5.15) and $0 \leq \mu < 1$
\begin{align*}
    A_2 (z,F) 
    &= \frac{\overline{\omega}e^{-\mu\pi i}}{(2\pi)^3 Qi} \sum_{k,n=1}^\infty \frac{{\mf}}{kn^2}  \times 2\pi i \left\{ (2\pi nkQ^2)^\frac{1}{2}\exp\left( \frac{\pi}{4}i + \frac{z}{2} + \frac{e^{-\frac{\pi}{2}i-z} }{4\pi nkQ^2} \right)  \right.
    \\
    &\times \left. \Gamma(2+\mu)\Gamma(2-\mu)W_{-\frac{3}{2}, \> \mu}\left( \frac{e^{-\frac{\pi}{2}i - z}}{2\pi nkQ^2} \right) - \sum_{m_1 = 0}^M R^{(4)}_{k,n, m_1}(z,\mu) - \sum_{m_2 = 0}^{M^\prime} R^{(4)}_{k,n, m_2}(z,\mu) \right\},
\end{align*}
\normalsize
where 
\begin{gather*}
    R^{(4)}_{k,n, m_1}(z,\mu)
    = \frac{(-1)^{m_1}}{m_{1}!}\Gamma(2\mu-m_1)\Gamma(2-\mu+m_1)(2\pi nkQ^2)^{\mu-m_1}e^{\left( z + \frac{\pi}{2}i \right)( \mu - m_1 )},  \tag{4.7.7}
    \\
    R^{(4)}_{k,n, m_2}(z,\mu)
    = \frac{(-1)^{m_2}}{m_{2}!}\Gamma(-2\mu-m_2)\Gamma(2+\mu+m_2)(2\pi nkQ^2)^{-\mu-m_2}e^{\left( z + \frac{\pi}{2}i \right)( -\mu - m_2 )},  \tag{4.7.8}
    \\
    R^{(4)}_{k,n,m,0}(z)
    =  \frac{m+1}{m!}\frac{e^{-\left( z+\frac{\pi}{2}i \right)m}}{(2\pi nkQ^2)^m}\left\{ \log(2\pi nkQ^2) + z  + \frac{\pi}{2}i + \sum_{k_1 = 1}^m \frac{1}{k_1} - C_0 - \frac{1}{m+1} \right\}  \tag{4.7.9}    
\end{gather*} and
\begin{align*}
   R^{(4)}_{k,n,m,\frac{1}{2}}(z) 
   &= \frac{\Gamma\left( \frac{3}{2}+m \right)}{(m!)^2}\frac{e^{\left( \frac{1}{2}-m \right)\left( z + \frac{\pi}{2} i \right)}}{(2\pi nkQ^2)^{m-\frac{1}{2}}}\left\{ m\left( \psi\left( \frac{3}{2}+m \right) -2\psi(m+1) \right)  \right.
   \\
   &- \left. m\left( \log (2\pi nkQ^2) + z + \frac{\pi}{2}i \right) + 1 \right\}   
     \tag{4.7.10}  
\end{align*}are residues of the integrand in $A_2 (z,F)$ at  $s = \mu - m_1, -\mu - m_2, -m$ and $\frac{1}{2}-m$ respectively. Finally, for $|y|<\pi, F \in \mathcal{S}^{\text{poly}}$ whose $(r,\lambda_j) = (1,1)$ in (1.5.15) and $0 \leq \mu < 1$, we have the series expression for $ f_1 (z,F) + f_{1}^{-}(z,F)$
\footnotesize
\begin{align*}
   f_1 (z,F) + f_{1}^{-}(z,F)
   &= \frac{\overline{\omega}e^{\mu\pi i}i}{(2\pi)^3 Q} \sum_{k,n=1}^\infty \frac{{\mf}}{kn^2} \times 2\pi i \left\{ (2\pi nkQ^2)^\frac{1}{2}\exp\left( -\frac{3}{4}\pi i + \frac{z}{2} + \frac{e^{\frac{3}{2}\pi i - z}}{4\pi nkQ^2} \right) \right. 
   \\
   &\times \left. \Gamma(2+\mu)\Gamma(2-\mu)W_{-\frac{3}{2}, \> \mu}\left( \frac{e^{\frac{3}{2}\pi i - z}}{2\pi nkQ^2} \right) - \sum_{m_1 = 0}^M R^{(1)}_{k,n, m_1}(z,\mu) - \sum_{m_2 = 0}^{M^\prime} R^{(1)}_{k,n, m_2}(z,\mu) \right\}
   \\
   &+\frac{\overline{\omega}e^{-\mu\pi i}}{(2\pi)^3 Qi} \sum_{k,n=1}^\infty \frac{{\mf}}{kn^2} \times 2\pi i \left\{ (2\pi nkQ^2)^\frac{1}{2}\exp\left( \frac{3}{4}\pi i + \frac{z}{2} + \frac{e^{-\frac{3}{2}\pi i - z} }{4\pi nkQ^2} \right)  \right.
   \\
    &\times \left. \Gamma(2+\mu)\Gamma(2-\mu)W_{-\frac{3}{2}, \> \mu}\left( \frac{e^{-\frac{3}{2}\pi i - z}}{2\pi nkQ^2} \right) - \sum_{m_1 = 0}^M R^{(2)}_{k,n, m_1}(z,\mu) - \sum_{m_2 = 0}^{M^\prime} R^{(2)}_{k,n, m_2}(z,\mu) \right\} 
    \\
    &-\frac{\overline{\omega}e^{\mu\pi i}}{(2\pi)^3 Qi} \sum_{k,n=1}^\infty \frac{{\mf}}{kn^2}  \times 2\pi i \left\{ (2\pi nkQ^2)^\frac{1}{2}\exp\left( -\frac{\pi}{4}i + \frac{z}{2} + \frac{e^{\frac{\pi}{2}i - z} }{4\pi nkQ^2} \right) \right. 
    \\
    &\times \left. \Gamma(2+\mu)\Gamma(2-\mu)W_{-\frac{3}{2}, \> \mu}\left( \frac{e^{\frac{\pi}{2}i - z}}{2\pi nkQ^2} \right) - \sum_{m_1 = 0}^M R^{(3)}_{k,n, m_1}(z,\mu) - \sum_{m_2 = 0}^{M^\prime} R^{(3)}_{k,n, m_2}(z,\mu) \right\}
    \\
    &+\frac{\overline{\omega}e^{-\mu\pi i}}{(2\pi)^3 Qi} \sum_{k,n=1}^\infty \frac{{\mf}}{kn^2}  \times 2\pi i \left\{ (2\pi nkQ^2)^\frac{1}{2} \exp\left( \frac{\pi}{4}i + \frac{z}{2} + \frac{e^{-\frac{\pi}{2}i - z} }{4\pi nkQ^2} \right) \right. 
    \\
    &\times \left. \Gamma(2+\mu)\Gamma(2-\mu)W_{-\frac{3}{2}, \> \mu}\left( \frac{e^{-\frac{\pi}{2}i - z}}{2\pi nkQ^2} \right) - \sum_{m_1 = 0}^M R^{(4)}_{k,n, m_1}(z,\mu) - \sum_{m_2 = 0}^{M^\prime} R^{(4)}_{k,n, m_2}(z,\mu) \right\}.   \tag{4.7.11}
\end{align*}
\normalsize
Since a lemma similar to Lemma 4.4.1 also holds, the third and the fourth series on the right hand side in (4.7.11) are absolutely and uniformly convergent on every compact subset on the whole complex plane. 

Next, by the theorem of residues, (2.2.7) and (4.6.5) we have
\begin{align*}
   f_2 (z,F) + f_{2}^{-}(z,F)
  &= \int_L \frac{\zeta(s-1)}{F(s)}e^{sz}ds - \int_{\overline{L}} \frac{\zeta(s-1)}{F(s)}e^{sz}ds
   \\
   &= -2\pi i \lim_{s \to 2} (s-2) \frac{\zeta(s-1)}{F(s)}e^{zs}
   \\
   &= -2\pi i\frac{e^{2z}}{F(2)}.  \tag{4.7.12}
\end{align*}Finally, by (4.3.1) and (4.6.6) 
\begin{equation*}
   f_3 (z,F) + f_{3}^{-}(z,F) = 0.  \tag{4.7.13}
\end{equation*}Thus, for $|y| < \pi$ we have
\begin{align*}
   2\pi i(f(z,F) + f^{-} (z,F))
   &= 2\pi i\left\{ \frac{1}{2\pi i}(f_1 (z,F) + f_{1}^{-}(z,F)) - \frac{e^{2z}}{F(2)} \right\}
   \\
   &= 2\pi iB(z,F),   \tag{4.7.14}
\end{align*}where
\begin{equation*}
   B(z,F) = \frac{1}{2\pi i}(f_1 (z,F) + f_{1}^{-}(z,F)) - \frac{e^{2z}}{F(2)}.  \tag{4.7.15}
\end{equation*}By (4.7.11), the function $f_1 (z,F) + f_{1}^- (z,F)$ is absolutely and uniformly convergent on every compact subset on the whole complex plane. Hence, the function $B(z,F)$ is an entire function. Since the function $f^- (z,F)$ has a meromorphic continuation for $y<\pi$ by Corollary 2.2.3, the function 
\[ f(z,F) = B(z,F) - f^- (z,F) \] is analytic for all $y<\pi$. Since the function $f(z,F)$ is analytic for $z \in \mathbb{H}$ and $y<\pi$, $f(z,F)$ can be analytically continued to the whole complex plane. In a similar manner, $f^- (z,F)$ can be analytically continued to the whole complex plane. Therefore, for all $z \in \mathbb{C}$ we have
\begin{equation*} 
   f(z,F) + f^{-}(z,F) =  B(z,F).  \tag{4.7.16}
\end{equation*}

Finally, we prove the functional equation (2.2.10). We recall the hypothesis that the coefficient $a_F (n)$  in the Dirichlet series of $F$ is real for all $n$. Hence, if $\rho$ is a non-trivial zero of $F$, then so is $\overline{\rho}$. For $z \in \mathbb{H}$ we have
\begin{align*}
   \overline{f(z,F)}
   &= \lim_{n \to \infty} \overline{\sum_{{ \scriptstyle \rho } \atop{\scriptstyle 0 < \text{Im} \hspace{0.01cm} \rho < T_n}} \frac{\zeta(\rho-1)}{F^\prime (\rho)}e^{\rho z}}
   \\
   &= \lim_{n \to \infty} \overline{\sum_{{ \scriptstyle \rho } \atop{\scriptstyle 0 < \text{Im} \hspace{0.01cm} \rho < T_n}} \left(\lim_{s \to \rho}\frac{s - \rho}{F(s) - F(\rho)} \zeta(s-1)e^{zs} \right)}.
\end{align*}Since $a_F (n) \in \mathbb{R}$, so $\overline{F(s)} = F(\overline{s})$ holds. Using this, we have
\begin{align*}
   f^- (\overline{z},F)
   &= \lim_{n \to \infty} \sum_{{ \scriptstyle \rho } \atop{\scriptstyle -T_n < \text{Im} \hspace{0.01cm} \rho < 0}} \frac{\zeta(\rho-1)}{F^\prime (\rho)}e^{\rho \overline{z}}
   \\
   &=  \lim_{n \to \infty} \overline{\sum_{{ \scriptstyle \rho } \atop{\scriptstyle -T_n < \text{Im} \hspace{0.01cm} \rho < 0}} \left(\lim_{s \to \rho}\frac{\overline{s} - \overline{\rho}}{\overline{F(s)} - \overline{F(\rho})} \overline{\zeta(s-1)}e^{z\overline{s}}\right)}
   \\
   &=  \lim_{n \to \infty} \overline{\sum_{{ \scriptstyle \rho } \atop{\scriptstyle -T_n < \text{Im} \hspace{0.01cm} \rho < 0}} \left(\lim_{s \to \rho}\frac{\overline{s} - \overline{\rho}}{{F(\overline{s})} - {F(\overline{\rho}})} \zeta(\overline{s}-1)e^{z\overline{s}} \right)}
   \\
   &= \lim_{n \to \infty} \overline{\sum_{{ \scriptstyle \rho } \atop{\scriptstyle 0 < \text{Im} \hspace{0.01cm} \rho < T_n}} \left(\lim_{\overline{s} \to \overline{\rho}}\frac{ s - \rho}{F(s) - F(\rho)} \zeta(s-1)e^{zs} \right)}
   \\
   &=  \lim_{n \to \infty} \overline{\sum_{{ \scriptstyle \rho } \atop{\scriptstyle 0 < \text{Im} \hspace{0.01cm} \rho < T_n}} \left(\lim_{s \to \rho}\frac{s - \rho}{F(s) - F(\rho)} \zeta(s-1)e^{zs} \right)}
   \\
   &= \overline{f(z,F)}.
\end{align*}Therefore, we have for $z \in \mathbb{H}$
\begin{equation*}
   f(z,F) = \overline{f^- (\overline{z},F)}.  \tag{4.7.17}
\end{equation*}Using (4.7.17), we have from (4.7.16)
\begin{align*}
   \overline{B(\overline{z},F)}
   &= \overline{f(\overline{z},F) + f^- (\overline{z},F)}
   \\
   &= \overline{f(\overline{z},F)} + \overline{f^- (\overline{z},F)}
   \\
   &=  \overline{f(\overline{z},F)} + f(z,F).
\end{align*}Using (4.7.17) and (4.7.16) again, we have for $z \in \mathbb{H}$
\begin{equation*} 
   \overline{f(\overline{z},F)} + f(z,F) = f^- (z,F) + f(z,F) = B(z,F).   \tag{4.7.18}
\end{equation*}Since $f(z,F), f^- (z,F)$ and $B(z,F)$ are entire functions, and (4.7.18) holds for all $z \in \mathbb{H}$, (4.7.18) holds for all $z \in \mathbb{C}$ by the analytic continuation. Therefore, the functional equation (2.2.10) holds for all $z \in \mathbb{C}$ and we have Theorem 2.2.4. $\square$  

\textbf{Acknowledgements}

Firstly, the author expresses his sincere gratitude to Prof. Kohji Matsumoto for many valuable comments and suggestions. If he did not accept the author, the author did not step in the world of mathematical research. The author would not be able to write this paper.
Secondly, the author would like to express his gratitude to Dr. Yuta Suzuki and Dr. Wataru Takeda for giving me various advice in obtaining Theorem 2.1.1 in Section 2.1. In particular, the author would like to express his gratitude to Dr. Sh\={o}ta Inoue for giving his many valuable comments on the whole of this paper. Thirdly, the author would like to express his gratitude to Prof. Isao Kiuchi for introducing the reference~\cite{Kac and Wie}. Fourthly,  the author expresses his gratitude to laboratory members, in particular, to master's student Keita Nakai giving the comment on the path of integration $\mathscr{L}$ in Section 4.2. Finally, the author expresses his gratitude to his parents for providing him with unfailing support and continuous encouragement.

\bigskip
　　　　　　　　　　　　　　Hideto Iwata 
                                                  
　　　　　　　　　　　　　　Graduate School of Mathematics

　　　　　　　　　　　　　　 Nagoya University
                                                   
　　　　　　　　　　　　　　 Furocho, Chikusa-ku,
                                                   
　　　　　　　　　　　　　　 Nagoya, 464-8602, Japan.

　　　　　　　　　　　　　　 \small{e-mail:d18001q@math.nagoya-u.ac.jp}

\end{document}